\documentclass[11pt]{article}
\usepackage{makeidx,amssymb,amscd}
\usepackage{amsthm}
\usepackage{amsmath}
\usepackage{amsfonts}
\usepackage{graphicx}
\usepackage{mathrsfs}
\usepackage{hyperref}
\usepackage{enumerate}
\usepackage{color}
\usepackage{latexsym}
\usepackage{caption}
\usepackage{subcaption}
\usepackage{enumerate}

\begin{document}
\numberwithin{equation}{section}
\newtheorem{axiom}{Axiom}[subsection]
\newtheorem{theorem}[axiom]{Theorem}
\newtheorem{lemma}[axiom]{Lemma}
\newtheorem{conjecture}[axiom]{Conjecture}
\newtheorem{proposition}[axiom]{Proposition}
\newtheorem{definition}[axiom]{Definition}
\newtheorem{remark}[axiom]{Remark}
\newtheorem{statement}[axiom]{Statement}
\newtheorem{claim}[axiom]{Claim}
\newtheorem{comment}[axiom]{Comment}
\newtheorem{corollary}[axiom]{Corollary}
\newtheorem{example}[axiom]{Example}
\newtheorem{question}[axiom]{Question}

\newcommand{\beq}{\begin{equation}}
\newcommand{\ee}{\end{equation}}
\newcommand{\bd}{\begin{definition}}
\newcommand{\ed}{\end{definition}}
\newcommand{\bt}{\begin{theorem}}
\newcommand{\et}{\end{theorem}}
\newcommand{\ble}{\begin{lemma}}
\newcommand{\ele}{\end{lemma}}
\newcommand{\bp}{\begin{proposition}}
\newcommand{\epo}{\end{proposition}}
\newcommand{\brm}{\begin{remark}}
\newcommand{\erm}{\end{remark}}
\newcommand{\bc}{\begin{corollary}}
\newcommand{\ec}{\end{corollary}}
\title{Survey on entropy-type invariants of sub-exponential growth in dynamical systems}
\author{Adam Kanigowski\footnote{A. K. was partially supported by the NSF grant DMS-1956310} \and Anatole Katok\footnote{A. K. was partially supported by the NSF grant DMS-16-02409} \and Daren Wei\footnote{D. W. was partially supported by the NSF grant DMS-16-02409}}
\maketitle
\begin{abstract}
Measure-theoretic and topological entropy are classical invariants in the theory of dynamical systems. There are several recently developed entropy type invariants for systems of sub-exponential growth: sequence entropy, slow entropy, Kakutani invariants, scaled entropy, entropy dimensions and entropy convergence rate. They measure the complexity of zero entropy systems by different approaches. These new invariants and corresponding new theories have many applications and interesting properties. This survey paper gives a comprehensive exposition of the slow entropy theory and also discusses some related topics.
\end{abstract}
\tableofcontents

\section{Introduction}

\subsection{Short history of entropy-type invariants}
Classification is one of  the central problems in ergodic theory and dynamical systems. For a dynamical system, one of its most important and natural characteristics  is its complexity. Therefore it is natural to try to classify dynamical systems based on their complexity. The classical and natural way to measure  orbit complexity for a dynamical system, continuous or measurable, is to  look at the growth of the number of  distinguishable orbits segments as a function of time. Distinguishable is the key word here and, according to  the setting, it is given a specific meaning. Basically there are three  principal settings for the problem corresponding to the following classes of dynamical systems: \emph{topological dynamical systems}, \emph{symbolic dynamical systems} and \emph{measurable dynamical systems}.

\paragraph{Distinguishing the orbit segments and counting.}\label{para:DistinguishingOrbitSegments} Orbit segments are  distinguished with the use of  auxiliary tools which can be viewed as mathematical abstraction of the notion of observation. Those tools are accordingly  \emph{Bowen metric} in the topological case, \emph{value of the zero coordinate} in the symbolic case, \emph{entropy $H(\xi)$\footnote{See \eqref{eq:entropyOfPartition} for detailed definition.}  of a partition} in the measurable case. Notice that the Bowen metric at orbit length $T$ is continuous as a function of $T$ and hence to distinguish orbit segments one needs to put an arbitrary threshold; in the  remaining cases the tool is discrete, but in the last one is rather arbitrary. To refine the device in topological case  the threshold is set  at an arbitrary low fixed level before  counting, and in  measurable one  more and more refined  partitions are considered, so that  the $\sigma$-algebra $\frak B(\xi)$ becomes more and more representative. In the symbolic case in principle one  may take values of several successive coordinates rather than a single one, but this  does not change the outcome for any   counting procedure. Similarly, for  many ``good''  continuous and measurable systems for adequate counting  it is sufficient to fix any  threshold  below a certain  fixed number (e.g. the {\em expansivity constant}), or to consider a single partition (e.g. a {\em generator}).

For the \emph{counting} we measure the asymptotic growth of orbit segments with the time length.  As a natural  measure of complexity one uses \emph{$(n,\epsilon)-$spanning set} or \emph{$(n,\epsilon)-$separated set}\footnote{See Section \ref{sec:topologicalentropy} for a detailed discussion.} in the topological setting, \emph{the number of sequences of successive $n$ symbols which appear in elements of a closed invariant set} in the symbolic setting and \emph{entropy $H(\xi^T_{-n})$}\footnote{$\xi^T_{-n}=\vee_{i=0}^{n-1}T^{-i}(\xi)$.} of the joint partition in the measurable setting.

\paragraph{Topological and measure-theoretic entropy.} The exponential growth of orbit complexity is given by the \emph{topological entropy} and \emph{measure-theoretic entropy}. These isomorphism invariants have been extensively studied over the last sixty years to deal with systems possessing exponential orbit growth. The basic connection between the entropies is provided by the Variational Principle for entropy\footnote{See Section \ref{sec:varpri}}. Based on Shannon-McMillan-Breiman theorem, one can interpret the variational  principle  as a precise quantification of  the following assertion:
\medskip

\noindent ($\cal VP$){\em The global  exponential growth rate of orbit complexity can be arbitrarily well approximated
by  the  exponential rates of growth  for the  sets of  statistically typical  orbits  satisfying various uniform distribution laws.}
\medskip

\paragraph{Distance between orbit segments.} The method of distinguishing orbit segments described above assumes uniform accuracy of the measuring device and, in the case of partitions, ignores its discrete nature. If we measure  complexity using partitions in the topological setting,  measurements  for nearby points and close to the boundary of the partition  may give different results.  These considerations motivate an alternative approach of measuring  the distance  between orbit segments:
using the average in time rather than maximal deviation.  As  it turns out,  the resulting  notions  do not change in the case of exponential  orbit complexity,  but  make crucial  difference   in the subexponential  case.  The integral or {\em Hamming metric}\footnote{See Definition \ref{def:HammingMetric} for general case.} between the  orbit segments  for a map $f$ or a flow $\Phi$ in a metric space  is defined  accordingly  as
\begin{equation}\label{eqdefintmetricmap}
d^{H,f}_n(x,y)=\frac1n\sum_{i=0}^{n-1} d(f^i(x),f^i(y))
\end{equation}
and
\begin{equation}\label{eqdefintmetricflow}
d^{H,\Phi}_t(x,y)=\frac1t\int_0^t d(\Phi^s(x),\Phi^s(y))\,ds.
\end{equation}

In the  case of discrete measurements,  which we use for symbolic and measurable systems, one  simply assumes that distance between symbols or different elements of the partition $\xi$ is equal to one and uses  \eqref{eqdefintmetricmap}. In other words, the Hamming distance  between two  $n$-strings of symbols is equal  to the proportion of the number of  places where  the symbols in the strings are different. Let $\Omega_N^n$  be the set of $n$-strings of sequences of symbols for an alphabet with $N$ elements. We
use simplified notation for the Hamming metric in this case:
\begin{equation}\label{Hammimgsymbolic}
d_n^H(\omega,\omega')= \frac1n\sum_{i=0}^{n-1}(1-\delta_{\omega_i\omega'_i}).
\end{equation}

Both the original measuring scheme  and the average one can be put into an even more general context by considering  a ``coherent '' system of metrics $d_n$ for the  sets $\Omega_N^n$ of $n$-strings of sequences of symbols for an alphabet with $N$ elements.  Coherence  means convexity with respect to  concatenation of strings.
(The original scheme of course  corresponds to  the trivial choice: $d_n(\omega, \omega')=1$ for $\omega\neq\omega'$).

A broad class of metrics  can be obtained by  fixing a set of elementary operations, and defining the distance as  the minimal
number of elementary operations needed to transform one string into the other, divided by the length of the strings. For the Hamming  metric elementary operations are just replacements of one symbol by another  at the same place.
Another  important special case is that of {\em Kakutani metrics} $d_n^M$ or \emph{$\bar{f}$-metric}\footnote{See \eqref{eq:fbarContinuous} for a continuous version.}, where an elementary operation consists in crossing out a symbol and inserting another one {\em at any place}.  Basic properties of those metrics are discussed in \cite[Section 4]{KatokMonotone}. For the case of
Hamming metrics for ergodic measure-preserving transformations this is  worked out in detail in \cite[Section 1]{KatokPubl}. The treatment in the Kakutani metric  case is quite similar. In  both cases  the key observation is that one can calculate the  size of an $\epsilon$-ball   using Stirling Formula and   on the exponential in $n$ scale this size goes to zero with $\epsilon$.  This also makes symbolic case straightforward.

\paragraph{Measuring subexponential growth.} In the case of zero topological entropy measuring the growth rates for $(n,\epsilon)-$spanning sets or $(n,\epsilon)-$separated sets and the number of sequences of successive $n$ symbols which appear in elements of a closed invariant sets provides a natural approach to the study of orbit complexity. In the symbolic case at least it is a well-researched  subject  branching into combinatorics, algorithmic complexity and suchlike, see \cite{Ferenczi2}. We will discuss the topological  case later.

However in the  measurable case  a difficulty appears. It turns out that  for any aperiodic measure preserving transformation $T$ any  speed of growth  below linear in $n$ appears for the entropy  $H(\xi^T_{-n})$ for an appropriately chosen  finite partition  $\xi$, see e.g.  \cite{Blume1}, although this is not the earliest proof of this fact in the literature.
Several solutions  have been suggested  to overcome this problem.

Before outlining those  let us point out  that independently of their success one cannot hope for a counterpart of variational principle
with  the standard choice of metric in the topological category since
topological orbit complexity  may change  for different uniquely ergodic realizations (e.g irrational rotation and the Denjoy example, \cite[Example 2.6.9]{HK-survey}). Hence a modification of the approach  in the topological setting may also be in order.

Classical entropy theory does not provide any information for systems with slower orbit growth, i.e. dynamical systems with sub-exponential orbit growth, for which the invariant is zero. In order to get some non-trivial invariants for such systems, one needs to observe their dynamical properties at different rates.

The earliest successful construction of slow growth  invariants  for measure preserving transformations was found by Kushnirenko \cite{Kushnirenko} and he denoted this invariant as the \emph{sequence entropy}. The idea here is to give the system enough time to develop exponential behavior, not against the number of iterates but rather against the number of measurements i.e.
instead of taking the standard consecutive sequence of iterates to compute the asymptotic growth of different codes, one calculates the asymptotic growth along a given (convenient) subsequence of consecutive iterates.

 Kushnirenko \cite{Kushnirenko} used the sequence entropy to show that the horocycle flow is not isomorphic with its cartesian square (see also \cite{KanigowskiVinhageWei}). Moreover, it is shown that this invariant vanishes at all scales if and only if the system has discrete spectrum. Following Kushnirenko's approach, Hulse \cite{Hulse}, Newton \cite{Newton1}, \cite{Newton2}, \cite{Newton3}, Krug \cite{KrugNewton}, Dekking \cite{Dekking} and Lema\'{n}czyk \cite{Lemanczyk} further developed the measure-theoretic sequence entropy and Goodman \cite{Goodman} modified this construction to the topological category. Sequence entropy is also related to classical ergodic theory properties, such as mixing : Saleski \cite{Saleski}, Hulse \cite{Hulse}, Huang, Shao and Ye \cite{HuangSahoYe};  mild mixing: Zhang \cite{Zhang}, Huang, Shao and Ye \cite{HuangSahoYe} and weak mixing: Hulse \cite{Hulse}, Zhang \cite{Zhang}, Huang, Shao and Ye \cite{HuangSahoYe}. The crucial idea behind the sequence entropy is that the lack of exponential growth is compensated with the sparseness of the defining subsequence.

On the other hand, in \cite{KatokTimeChange}, \cite{KatokMonotone}  A. Katok suggested to count the number of statistically
different orbits (Hamming balls) in a scale which can be adjusted to the system (subexponential,
polynomial, logarithmic, etc.). This approach was further developed by Katok and Thouvenot \cite{KatokThou}, where the classical slow entropy definition was generalized to actions of amenable groups and then was used to give a criterion for smooth realization of $\mathbb{Z}^k$ actions. Then Galatolo \cite{Galatolo} generalized this invariant to the topological setting. Later Ferenczi \cite{Ferenczi} proved that vanishing of the measure-theoretic slow entropy at all scales is equivalent to discrete spectrum in measurable category; recently, Kanigowski, Vinhage and Wei \cite{KanigowskiVinhageWei} proved an analogue of this result in the topological category. Hochman \cite{Hochman} studied slow entropy for higher rank actions and gave a counterexample for smooth realization in the setting of infinite-measure preserving $\mathbb{Z}^2$-actions. A.\ Katok, S. Katok and Rodriguez-Hertz \cite{KatokKatokRH} studied the connections between the slow entropy and Fried entropy for maximal rank $\mathbb{Z}^k$-actions as well as the lower bounds of these quantities for maximal rank $\mathbb{Z}^k$-actions.

One of the most important differences between slow entropy and sequence entropy is that slow entropy is better adjusted to deal with non-homogeneous systems. For homogeneous systems, computation of both invariants is based on controlled (polynomial) divergence of nearby points (see \cite{KanigowskiVinhageWei} for more details). Therefore, it can be shown that the metric orbit growth is the same (up to constants) as the topological orbit growth and the second one is easier to compute while dealing with, say, unipotent matrices. However, if the system lacks homogeneity, i.e. some rank one systems, translation flows, smooth flows on surfaces, Liouvillian time changes of linear flows on $\mathbb{T}^2$, time changes of unipotent flows, then the sequence entropy is very difficult to control. Indeed, since one only takes a (sparse) subsequence of iterations, it is almost impossible to control what happens between two consecutive elements of the subsequence. On the other hand, slow entropy tracks entire pieces of orbits and therefore it is easier to control the growth.

The aim of this survey is to give a systematic exposition of slow entropy type invariants. Most of the survey is devoted to the slow entropy introduced by Katok and Thouvenot. Besides slow entropy and sequence entropy, we will also describe several other invariants which share some common features with slow entropy, i.e. Fried entropy (\cite{KatokKatokRH}), Kakutani equivalence invariants (\cite{RatnerKakInv}, \cite{RatnerKakInv2}, \cite{KanigowskiVinhageWei2}), scaled entropy (\cite{ZhaoPesin}, \cite{ZhaoPesin2}), entropy dimensions (\cite{DeCarvalho}, \cite{FerencziPark}, \cite{ADP}, \cite{DouHuangPark}, \cite{DouHuangPark2}) and entropy convergence rates (\cite{Blume1}, \cite{Blume2}, \cite{Blume3}). These invariants reflect different features of a given dynamical system from different points of view.

\subsection{Plan of the paper}
The structure of this survey is as follows:
\begin{itemize}
\item In $\mathsection\ref{sec:BTTM}$ we review basic theory of entropy, which includes the topological entropy, measure-theoretic entropy and variational principle.
\item In $\mathsection\ref{sec:SequenceEntropy}$ we introduce the sequence entropy, describe its properties and connections with the classical entropy. Moreover, we mention spectral theory, variational principle and provide some examples of the sequence entropy in homogeneous setting.
\item In $\mathsection\ref{sec:slowMetricEntropy}$ we give a systematic introduction of the slow entropy, relations with spectral theory, Kushnirenko's inequality, conditions on vanishing of slow entropy, failure of variational principle and different upper and lower quantities of slow entropy.. Several interesting examples are also provided to illustrate the features of slow entropy for example, quasi-unipotent flows, rank one systems, surface flows, AbC constructions and some $\mathbb{Z}^k-$actions.
\item In $\mathsection\ref{sec:Otherinvariants}$ we mainly concentrate on several invariants related to slow entropy, i.e. Fried entropy, Kakutani equivalence invariants, scaled entropy, entropy dimensions and entropy convergence rate. Several examples of these invariants with detailed discussions are also provided in this section.
\item In $\mathsection\ref{sec:openQuestions}$ we post some open questions related to slow entropy invariants.
\end{itemize}

\textbf{Acknowledgements.} The authors would like to thank Svetlana Katok and Mariusz Lema\'{n}czyk for their help, support and careful reading of the first draft of the paper.

\section{Classical entropy theory}\label{sec:BTTM}
In this section, we will review the classical theory of topological and measure-theoretic entropy, which will help us formulate parallel results for the sequence entropy and the slow entropy. For more details see e.g.\ \cite{KatokHassel}.
\subsection{Topological entropy}\label{sec:topologicalentropy}
Let $(X,d)$ be a compact metric space and $T:(X,d)\to (X,d)$ be a continuous map. For $\epsilon>0$ let $S_d(\epsilon)$  denote  the minimal number of balls of radius $\epsilon$ which covers the whole space, $D_d(\epsilon)$  the minimal number of sets with diameter less or equal than $\epsilon$ whose union covers the whole space and $N_d(\epsilon)$ the maximal number of $\epsilon$-separated points in the space. It is clear that these three quantities satisfy the following inequalities:
\begin{equation}\label{eq:topologicalBallIneq}
\begin{aligned}
D_d(2\epsilon)\leq S_d(\epsilon)\leq D_d(\epsilon),\\
N_d(2\epsilon)\leq S_d(\epsilon)\leq N_d(\epsilon).
\end{aligned}
\end{equation}

In order to measure the complexity of the dynamical system $(X,T,d)$, the Bowen metric is introduced as follows:
\begin{equation}\label{eq:BowenMetric}
d_n^T(x,y):=\max\{d(x,y),d(Tx,Ty),\ldots,d(T^{n-1}x,T^{n-1}y)\}.
\end{equation}

Combining the quantities defined in the first paragraph and the Bowen metric, the complexity of $(X,T,d)$ with respect to $\epsilon>0$ is measured by the following quantity:
\begin{equation}
h_d(T,\epsilon):=\limsup_{n\to\infty}\frac{\log S_{d_n^T}(\epsilon)}{n}.
\end{equation}

Notice that $h_d(T,\epsilon)$ is non decreasing as a function of $\epsilon$ and thus we can define the quantity $h_d(T)$ of $(X,T,d)$ as,
\beq
h_d(T):=\lim_{\epsilon\to0}h_d(T,\epsilon).
\ee

In fact, $h_d(T)$ is independent of the  metric $d$:
\bp\label{prop:topologicalEntropyIndemetric}
If $d'$ is another metric on $X$ which defines the same topology as $d$, then $h_d(T)=h_{d'}(T)$.
\epo

Thus we have following,
\bd[Topological entropy]
The quantity $h_d(T)$ calculated for any metric generating a given topology in $X$ is called the topological entropy of $T$ and is denoted as $h_{top}(T)$.
\ed
\brm
We point out that the following quantities\footnote{The reason they are equivalent are inequalities \eqref{eq:topologicalBallIneq}.} are equivalent to the definition of topological entropy,
\begin{equation}
\begin{aligned}
h_{top}(T)&=\lim_{\epsilon\to0}\liminf_{n\to\infty}\frac{\log S_{d_n^T}(\epsilon)}{n}=\lim_{\epsilon\to0}\lim_{n\to\infty}\frac{\log D_{d_n^T}(\epsilon)}{n}\\&=\lim_{\epsilon\to0}\limsup_{n\to\infty}\frac{\log N_{d_n^T}(\epsilon)}{n}=\lim_{\epsilon\to0}\liminf_{n\to\infty}\frac{\log N_{d_n^T}(\epsilon)}{n}.\\
\end{aligned}
\end{equation}
\erm

The following proposition describes some basic features of topological entropy:
\begin{proposition}\label{prop:TopologicalEntropyProp}
We have the following basic properties of topological entropy:

(1) If the map $S$ is a factor of $T$, then $h_{top}(S)\leq h_{top}(T)$.

(2) If $\Lambda$ is a closed $f-$invariant set, then $h_{top}(T{\upharpoonright\Lambda})\leq h_{top}(T)$.

(3) If $X=\cup_{i=1}^m\Lambda_i$, where $\Lambda_i$, $i=1,2,\ldots,m$ are closed $T-$invariant sets, then $h_{top}(T)=\max_{1\leq i\leq m}h_{top}(T_{\upharpoonright\Lambda_i}).$

(4) $h_{top}(T^m)=|m|h_{top}(T)$.

(5) $h_{top}(T\times S)=h_{top}(T)+h_{top}(S)$.
\end{proposition}

We end this section with the following basic examples.
\begin{example}
1. The topological entropy of any translation $T_{\gamma}$ of the torus or any linear flow $(T_w^t)$ on the torus is equal to zero.

2. The topological entropy of the gradient flow on the round sphere is equal to zero.

3. If $E_m:S^1\to S^1$ is an expanding map of degree $m$, $E_m(z)=z^m$, then
$h_{top}(E_m)=\log|m|$.

4. For a topological Markov chain $\sigma_A$, $h_{top}(\sigma_A)=\log|\lambda_A^{\max}|$, where $\lambda_A^{\max}$ is the maximal eigenvalue of $A$.
\end{example}
\subsection{Measure-theoretic entropy}
Instead of considering the complexity in the topological category, similar invariant can be constructed in the measurable category for measure preserving  dynamical systems.

Let $(X,\mathscr{B},\mu)$ be a probability Borel space and $I$ a finite or countable set of indices. Suppose we have a measurable partition $\xi=\{C_{\alpha}|\alpha\in I\}$. We define $H(\xi)$ as:

\begin{equation}\label{eq:entropyOfPartition}
H(\xi)=H_{\mu}(\xi)=-\sum_{\alpha\in I,\mu(C_{\alpha}>0)}\mu(C_{\alpha})\log\mu(C_{\alpha}).
\end{equation}
The \emph{conditional entropy} of $\xi$ with respect to a partition $\eta=\{D_{\alpha}|\alpha\in J\}$ is defined as
$$H(\xi|\eta)=-\sum_{\beta\in J}\mu(D_{\beta})\sum_{\alpha\in I}\mu(C_{\alpha}|D_{\beta})\log(\mu(C_{\alpha}|D_{\beta})),$$
where $\mu(A|B)=\frac{\mu(A\cap B)}{\mu(B)}$.

Moreover, we introduce the following notions for partitions:
\begin{enumerate}[(1)]
  \item $\xi\leq\eta$ if and only if for all $D\in\eta$ there exists a $C\in \xi$ such that $D\subset C$ and we will say $\eta$ is a \emph{refinement} of $\xi$ and $\xi$ is subordinate to $\eta$;
  \item The \emph{joint partition} of $\xi$ and $\eta$ is defined as $\xi\vee\eta:=\{C\cap D|C\in\xi,\ \ D\in\eta,\ \ \mu(C\cap D)>0\};$
  \item Two partitions $\xi$ and $\eta$ are \emph{independent} if $\mu(C\cap D)=\mu(C)\cdot\mu(D)$ for all $C\in\xi$, $D\in\eta$.
\end{enumerate}

We have the following proposition.
\bp\label{prop:basicEntropy} Let $\xi=\{C_{\alpha}|\alpha\in I\}$, $\eta=\{E_{\alpha}|\alpha\in J\}$, $\zeta=\{E_{\alpha}|\alpha\in K\}$ be finite or countable measurable partitions of $(X,\mu)$. Then:
\begin{itemize}
  \item[(1)] $0<-\log(\sup_{\alpha\in I}\mu(C_{\alpha}))\leq H(\xi)\leq\log\operatorname{card}\xi$; furthermore if $\xi$ is finite then $H(\xi)=\log\operatorname{card}\xi$ if and only if all elements of $\xi$ have equal measure.
  \item[(2)] $0\leq H(\xi|\eta)\leq H(\xi)$; $H(\xi|\eta)=H(\xi)$ if and only if $\xi$ and $\eta$ are independent; $H(\xi|\eta)=0$ if and only if $\xi\leq\eta(\textrm{mod}0)$. If $\zeta\geq\eta$ then $H(\xi|\zeta)\leq H(\xi|\eta)$.
  \item[(3)] $H(\xi\vee\eta|\zeta)=H(\xi|\zeta)+H(\eta|\xi\vee\zeta)$. In particular, for $\zeta=\nu$ we obtain
      \beq\label{conditionalEntropy}
      H(\xi\vee\eta)=H(\xi)+H(\eta|\xi).
      \ee
  \item[(4)] $H(\xi\vee\eta|\zeta)\leq H(\xi|\zeta)+H(\eta|\zeta)$; in particular $H(\xi\vee\eta)\leq H(\xi)+H(\eta)$.
  \item[(5)] $H(\xi|\eta)+H(\eta|\zeta)\geq H(\xi|\zeta)$.
  \item[(6)] If $\lambda$ is another measure on $X$ then for every partition $\xi$ measurable for both $\mu$ and $\lambda$ and for any $p\in[0,1]$,
      $$pH_{\mu}(\xi)+(1-p)H_{\lambda}(\xi)\leq H_{p\mu+(1-p)\lambda}(\xi).$$
\end{itemize}
\epo

Assume now we have a measure-preserving dynamical system $(X,\mathscr{B},\mu,T)$. The definition of $H(\xi)$ and Proposition \ref{prop:basicEntropy} allow us to measure the complexity of $(X,\mathscr{B},\mu,T)$ with respect to a partition $\xi$. Define the $\xi^T_{-n}=\bigvee_{i=0}^{n-1}T^{-i}(\xi)$, then $H_n(T,\xi)$ as follows:
\beq
H_n(T,\xi)=-\frac{1}{n}\sum_{C\in\xi^T_{-n}}\mu(C)\log\mu(C).
\ee

By subadditivity property of $H(\xi)$ (Proposition \ref{prop:basicEntropy} $(4)$), we have $$(n+m)H_{n+m}(T,\xi)\leq nH_{n}(T,\xi)+mH_m(T,\xi),$$ and thus
$$h(T,\xi):=\lim_{n\to\infty}H_n(T,\xi)
$$
exists.
\brm
$h(T,\xi)$ is also equal to the limit of the conditional entropy:
\beq
h(T,\xi)=\lim_{n\to\infty}H(\xi|T^{-1}(\xi_{-n}^T)).
\ee
\erm

In order to simplify our notation, we define the \textit{Rokhlin distance} between two partitions $\xi$ and $\eta$ as
\begin{equation}\label{eq:RokhlinDistance}
\rho(\xi,\eta)=H(\xi|\eta)+H(\eta|\xi).
\end{equation}

\bp\label{prop:MetricEntropyProp1}
The following are some basic properties of $h(T,\xi)$:
\begin{itemize}
\item[(1)] $0\leq\limsup_{n\to\infty}-\frac{1}{n}\log(\sup_{C\in\xi^T_{-n}}\mu(C))\leq h(T,\xi)\leq H(\xi)$.

\item[(2)] $h(T,\xi\vee\eta)\leq h(T,\xi)+h(T,\eta)$.

\item[(3)] $h(T,\eta)\leq h(T,\xi)+H(\eta|\xi)$; in particular if $\xi\leq\eta$ then $h(T,\xi)\leq h(T,\eta)$.

\item[(4)] $|h(T,\xi)-h(T,\eta)|\leq \rho(\xi,\eta)$ (the Rokhlin inequality).

\item[(5)] $h(T,T^{-1}(\xi))=h(T,\xi)$ and if $T$ is invertible $h(T,\xi)=h(T,T(\xi))$.

\item[(6)] $h(T,\xi)=h(T,\vee_{i=0}^kT^{-i}(\xi))$ for $k\in \mathbb{N}$ and if $T$ is invertible $h(T,\xi)=h(T,\vee_{i=-k}^kT^i(\xi))$ for all $k\in \mathbb{N}$.
\end{itemize}
\epo


By taking supremum of $h(T,\xi)$ over all measurable finite entropy partitions, the entropy of $(X,\mathscr{B},\mu,T)$ is defined as follows:
\bd[Measure-theoretic entropy]
The measure-theoretic entropy of $T$ with respect to $\mu$ is
\begin{equation}
h_{\mu}(T):=\sup\{h(T,\xi)\;|\;\xi\text{ is a measurable partition with }H(\xi)<\infty\}.
\end{equation}
\ed

One of the most important features of measure-theoretic entropy is the generator theorem, which is a crucial property for applications.

\bd
A partition $\xi$ is called a generator for $T$ if $\mathfrak{E}=\{\xi\}$ is a sufficient family. A family $\mathfrak{E}$ of measurable partitions with finite entropy is called sufficient with respect to the measure-preserving transformation $T$ if

(1) for a noninvertible $T$, partitions subordinate to partitions of the form $\vee_{i=0}^kT^{-i}(\xi)$ ($\xi\in\mathfrak{E}$, $k\in \mathbb{N}$) form a dense subset in the space of all partitions with finite entropy equipped with the Rokhlin metric.

(2) for an invertible $T$ the same holds for partitions subordinate to $\vee_{i=-l}^lT^i(\xi)$ ($\xi\in\mathfrak{E}$, $l\in \mathbb{N}$).
\ed

\bt\label{GeneratorTheorem}
If $\xi$ is a generator for $T$ then $h_{\mu}(T)=h_{\mu}(T,\xi)$. More generally, $h_{\mu}(T)=\sup_{\xi\in\mathfrak{E}}h_{\mu}(T,\xi)$ for any sufficient family $\mathfrak{E}$ of partitions.
\et

It is worth to notice that the measure-theoretic entropy can be introduced in a similar way as topological entropy. More precisely, the measure-theoretic entropy turns out to be the asymptotic value of the number of balls needed to cover a subset of positive measure instead of the whole space $X$ \cite{KatokPubl}: let $N(T,n,\epsilon,\delta)$ denote the minimal number of $\epsilon-$balls in the $d_n^T-$metric which cover a set of measure at least $1-\delta$.

\bt[Katok \cite{KatokPubl}]\label{MetricEntropyBall}
For every $\delta>0$,
$$h_{\mu}(T)=\lim_{\epsilon\to0}\liminf_{n\to\infty}\frac{\log N(T,n,\epsilon,\delta)}{n}=\lim_{\epsilon\to0}\limsup_{n\to\infty}\frac{\log N(T,n,\epsilon,\delta)}{n}.$$
\et


We list some basic properties of measure-theoretic entropy below:
\begin{proposition}\label{prop:MetricEntropyProp}
Let $T:(X,\mu)\to(X,\mu)$ be a measure-preserving transformation of a probability space $(X,\mathscr{B},\mu)$ and $\eta,\xi$ be measurable partitions with finite entropy. Then:
\begin{itemize}
\item[(1)] If $S:(Y,\nu)\to(Y,\nu)$ is a factor of $T:(X,\mu)\to(X,\mu)$ then $h_{\mu}(S)\leq h_{\mu}(T)$.

\item[(2)] If $A$ is invariant for $T$ and $\mu(A)>0$ then $h_{\mu}(T)=\mu(A)h_{\mu_A}(T)+\mu(X\backslash A)h_{\mu_{X\backslash A}}(T)$.

\item[(3)] If $\mu,\lambda$ are two invariant probability measures for $T$ then for any $p\in[0,1]$
$$h_{p\mu+(1-p)\lambda}(T)\geq ph_{\mu}(T)+(1-p)h_{\lambda}(T).$$

\item[(4)] $h_{\mu}(T^k)=kh_{\mu}(T)$ for any $k\in \mathbb{N}$. If $T$ is invertible then $h_{\mu}(T^{-1})=h_{\mu}(T)$ and hence $h_{\mu}(T^k)=|k|h_{\mu}(T)$ for any $k\in \mathbb{Z}$.

\item[(5)] $h_{\mu\times\lambda}(T\times S)=h_{\mu}(T)+h_{\lambda}(S)$.
\end{itemize}
\end{proposition}

We recall the following simple examples to end this section.
\begin{example}
1. Rotations and linear flows on tori have zero measure-theoretic entropy (with respect to Lebesgue measure $\lambda$)\footnote{The entropy of a flow is defined as the entropy of its time $1$ map.}.

2. For an expanding map $E_k$ on $S^1$, its measure-theoretic entropy (with respect to Lebesgue measure $\lambda$) is
$$h_{\lambda}(E_k)=\log|k|.$$

3. For the full $N$-shift $\sigma_N$ with Bernoulli measure $(p_1,\ldots,p_N)$, its measure-theoretic entropy is
$$h_{\mu}(\sigma_N)=-\sum_{i=1}^{N}p_i\log p_i.$$
\end{example}

\subsection{Variational Principle}\label{sec:varpri}
Once we have topological entropy and measure-theoretic entropy, it is natural to ask whether there is any relation between these two invariants. Notice first that there are some important differences between these two quantities: the measure-theoretic entropy of the union of two invariants sets is the weighted sums of measure-theoretic entropies but the topological entropy is the maximum of the two topological entropies. This corresponds to the observation that topological entropy measures the maximal complexity and measure-theoretic entropy measures the statistical complexity (statistics based on a given measure). The classical relation between these two entropies is given by the variational principle:


\bt[Variational Principle]
If $f:X\to X$ is homeomorphism of a compact metric space $(X,d)$, then
\beq
h_{top}(f)=\sup \{h_{\mu}(f)|\mu\in\mathfrak{M}(f)\}.
\ee
Here $\mathfrak{M}(f)$ is the set of all $f-$invariant Borel probability measures.
\et



It is natural to ask whether there always exists a measure that achieves the above supremum. In general, this is not true see e.g.\ \cite{Buzzi}, \cite{Misiu} and \cite{Ruette}. Below, we make some comments on the existence and uniqueness of measures of maximal entropy.

\brm
Expansive maps of compact metric spaces have a measure of maximal entropy.
\erm

\brm
Measure of maximal entropy might not be unique. The simplest example is to consider a union of several disjoint copies of the same expansive system (which is obviously expansive). In this case, there are more than one measure of maximal entropy.
\erm
\brm
All transitive topological Markov chains, horseshoes, hyperbolic toral automorphisms have a unique measure of maximal entropy.
\erm

We provide several examples at the end of this section.

\begin{example}
1. By Variational Principle, it is easy to see the following formula,
$$h_{top}(f)=h_{top}(f_{\upharpoonright_{NW(f)}}).$$
where $NW(f)$ is the set of all nonwandering points of $f$, where we recall that  $x\in X$ is a nonwandering point of $f$ if for any open set $U\i x$ there exists an $N>0$ such that $f^N(U)\cap U\neq\emptyset$. The reason we have the above equation is that the support of any invariant measure is included in $NW(f)$ and we then use Variational Principle.

2. By the previous example it follows that the measure-theoretic entropy of the gradient flow for any function with isolated critical points is zero.

3. By Variational Principle, the maximal entropy of a full shift (with $N$ symbols) is $\log N$ and the maximal entropy measure is $(\frac{1}{N},\ldots,\frac{1}{N})$.
\end{example}

\subsection{Kushnirenko's inequality}\label{sec:KushnirenkoSection}
In the theory of dynamical systems, one of the most challenging open problems is that of smooth realization of measure-preserving transformations. The only known restriction for the  smooth realization problem is the finite entropy condition, which was originally introduced by Kushnirenko for absolutely continuous measures. We will now recall Kushnirenko's inequality.
We need to recall several definitions:
\bd
Let $(X,d)$ be a compact metric space and $b(\epsilon)$ the minimal cardinality of a covering of $X$ by $\epsilon-$balls. Then
$$D(X)=\limsup_{\epsilon\to0}\frac{\log b(\epsilon)}{|\log\epsilon|}\in \mathbb{R}\cup\{\infty\}$$
is called the ball dimension of $X$.
\ed
\bd
Let $(X,d)$ be a compact metric space and $f:X\to X$ a Lipschitz  map. Then the Lipschitz constant $L(f)$ of $f$ is defined by
$$L(f)=\sup_{x\neq y}\frac{d(f(x),f(y))}{d(x,y)}.$$
\ed

\bt[Kushnirenko's inequality]\label{thm:KushnirenkoInequality}
Let $(X,d)$ be a compact metric space of finite ball dimension $D(X)$ and $f:X\to X$ be a Lipschitz map. Then
\beq
h_{top}(f)\leq D(X)\max(0,\log L(f)).
\ee
\et

A natural corollary of Kushnirenko's inequality together with the variational principle (see Section \ref{sec:varpri}) is the following obstruction for smooth realization:
\bc
For a $C^1$ map $f:M\to M$ of a compact Riemannian manifold
$$h_{top}(f)\leq\max(0,\dim M\log\sup_x\|Df_x\|)<\infty.$$
\ec

\subsection{Mixing, weak mixing and mild mixing}
Sequence entropy and slow entropy are related to some ergodic and spectral properties. We will introduce basic notions both in the measurable and topological category.

\bd[Mixing, weak mixing, rigid and mild mixing]\label{def:MixingRigid}
Let  $(X,\mathscr{B},\mu,T)$  be a measure preserving dynamical system. Then
\begin{enumerate}[(a)]
\item Let $U_T:L^2(X,\mu)\to L^2(X,\mu)$ denote the \emph{Koopman operator} associated with $T$, i.e. $U_T(f)=f\circ T$;
\item $T$ is \emph{mixing} if for all $B,C\in\mathscr{B}$:
\beq
\lim_{n\to\infty} \mu(T^{-n}B\cap C)=\mu(B)\mu(C);
\ee
\item $T$ is \emph{weakly mixing} (or has \emph{continuous spectrum}) if for all $B,C\in\mathscr{B}$:
\beq
\lim_{n\to\infty}\frac{1}{n}\sum_{i=0}^{n-1}|\mu(T^{-i}B\cap C)-\mu(B)\mu(C)|=0;
\ee
\item $T$ is \emph{rigid} if there exists an increasing sequence $\{t_n\}$ such that for any $f\in L^2(X,\mathscr{B},\mu)$, $$T^{t_n}f\to f$$ in $L^2$;

\item $T$ is \emph{mildly mixing} if there is no non-constant rigid functions in $L^2(X,\mathscr{B},\mu)$, where a function $f\in L^2(X,\mathscr{B},\mu)$ is called rigid if there exists an increasing sequence $\{t_n\}$ such that $T^{t_n}f\to f$ in the $L^2-$topology.

\end{enumerate}
\ed
We have the following inclusions:
\beq\label{eq:MixingRelations}
\textrm{Mixing systems}\subset\textrm{Mildly mixing systems}\subset\textrm{Weakly mixing systems}\subset\textrm{Ergodic systems}.
\ee
The first inclusion is straightforward from the above definitions. The second inclusion is based on the following characterizations of weak mixing and mild mixing:
\bt[Furstenberg-Weiss \cite{FurstenburgWeiss}]
$T$ is mildly mixing if for all ergodic systems $(Y,S,\nu,\mathscr{C})$ (where $\nu(Y)$ might be infinite) not of type $I$, the product system $(X\times Y, \mathscr{B}\times\mathscr{C}, \mu\times\nu,T\times S)$ is also ergodic.
\et
We recall that $T$ is weakly mixing if and only if for all \textbf{finite} ergodic systems $(Y,S,\nu,\mathscr{C})$, the product system $(X\times Y, T\times S, \mu\times\nu, \mathscr{B}\times\mathscr{C})$ is also ergodic. Moreover, we also recall that all inclusions in \eqref{eq:MixingRelations} are strict.

The above definitions have their analogs in the topological category:
\bd[Topological properties: transitivity, mixing, weak mixing and mild mixing]\label{def:topologicalSettingMixing}
Let $(X,T)$ be a topological dynamical system:
\begin{enumerate}[(a)]
\item $(X,T)$ is said to be \emph{transitive} if any closed invariant subset is either nowhere dense or the whole space $X$;
\item $(X,T)$ is said to be \emph{topologically weakly mixing} if $(X\times X, T\times T)$ is transitive;
\item $(X,T)$ is said to be \emph{topologically mixing} if for any non-empty open subsets $U$ and $V$ of $X$, there is $N\in\mathbb{N}$ such that $U\cap T^{-n}V\neq\emptyset$ for each $n\geq N$;
\item $(X,T)$ is said to be \emph{topologically mildly mixing} if for any transitive system $(Y,S)$, $(X\times Y, T\times S)$ is transitive.
\end{enumerate}
\ed


\section{Sequence entropy}\label{sec:SequenceEntropy}
In this section, we introduce the notion of sequence entropy and describe its basic properties. We also discuss relations between sequence entropy and classical entropy. Several examples are provided to illustrate some interesting properties of sequence entropy.

\subsection{Measure-theoretic sequence entropy}
Notice that a crucial feature of classical entropy is that it is an  isomorphism invariant\footnote{Topological entropy is an isomorphic invariant due to Proposition \ref{prop:TopologicalEntropyProp} $(1)$ and measure-theoretic entropy is a metric invariant due to Proposition \ref{prop:MetricEntropyProp} $(1)$.}. As a consequence,  entropy is a powerful tool in the (isomorphism) classification problem of dynamical systems. For example, by calculating measure-theoretic entropy, we can easily show that Bernoulli shifts $(\frac{1}{2},\frac{1}{2})$ and $(\frac{1}{3}, \frac{1}{3}, \frac{1}{3})$ are not isomorphic. However, for systems of zero entropy classical entropy theory provides no information on classification\footnote{For example the horocycle flow and its Cartesian square.}. In order to distinguish non-isomorphic  zero entropy systems, Kushnirenko \cite{Kushnirenko} developed an invariant based on modifying the classical entropy definition by replacing iterates along the sequence of natural numbers by iterates along a fixed subsequence.
\begin{definition}[Sequence entropy, Kushnirenko \cite{Kushnirenko}]\label{def:sequenceMetricEntropy}
Suppose that $(X,\mathscr{B},\mu,T)$ is a dynamical system (where $T$ is invertible). For a given integer-valued sequence $A=\{t_1,t_2,\ldots,t_n,\ldots\}$, a measurable partition (finite or countable) $\xi$ satisfying $H(\xi)<\infty$, let
\begin{equation}
\begin{aligned}
&h_A(T,\xi):={\limsup_{n\to\infty}}\frac{1}{n}H(\bigvee_{i=1}^nT^{-t_i}\xi)\\
&h_{A,\mu}(T):=\sup_{\xi}h_A(T,\xi).
\end{aligned}
\end{equation}

We call $h_{A,\mu}(T)$  the (measure-theoretic) sequence entropy of $T$.
\end{definition}

By an analogous argument as for the classical entropy, one shows that  $h_{A,\mu}(T)$ is an isomorphism invariant. If we take $t_i=i-1$, $h_{A,\mu}(T)$ coincides with the classical entropy. However, by adjusting the subsequence to a given dynamical system (eg. a system of zero entropy), sequence entropy can often provide additional information about the dynamics. We point out an interesting property of sequence entropy: Dekking \cite{Dekking} proved that measure-theoretic sequence entropy does not depend monotonically on $A$: if $A$ is a subsequence of $B$, it is possible that $h_{A,\mu}(T)<h_{B,\mu}(T)$.

By adjusting the proofs for classical entropy, one recovers analogous properties of the sequence entropy:

\begin{proposition}[Kushnirenko \cite{Kushnirenko}]\label{prop:RokhlinSequence}
\begin{equation}
|h_{A,\mu}(T,\xi)-h_{A,\mu}(T,\eta)|\leq\rho(\xi,\eta),
\end{equation}
where $\rho(\xi,\eta)$ is defined in \eqref{eq:RokhlinDistance}.
\end{proposition}

We also have a generator-type result:
\begin{proposition}[Kushnirenko \cite{Kushnirenko}]
Let $\xi_1\leq\xi_2\leq\ldots\leq\xi_n\leq\ldots$ and $\bigcap_{i=1}^{\infty}\xi_i=\epsilon$ ($\epsilon$ is the point partition). Then for any $A$ and $T$, we have,
\begin{equation}
h_{A,\mu}(T)=\lim_{k\to\infty}h_{A,\mu}(T,\xi_k).
\end{equation}
\end{proposition}

Once we can control the growth rate of the transformation $T$, we have:
\begin{proposition}[Kushnirenko \cite{Kushnirenko}]
Suppose that $T$ is a diffeomorphism of a compact $m-$dimensional Riemannian manifold $M^m$ and $\|T\|$ be the maximum expansion of an $(m-1)-$dimensional non-orientated volume along all $(m-1)-$dimensional tangential directions. Let  $A=\{t_1,t_2,\ldots,t_n,\ldots\}$ be a given sequence. If $\|T^{t_n}\|\leq C\lambda^n$ (where $C>0$ does not depend on $n$), then $h_{A,\mu}(T)\leq m\log\lambda$.
\end{proposition}

However, the sequence entropy behaves differently than the classical entropy when considering cartesian products:
\begin{proposition}[Kushnirenko \cite{Kushnirenko}, Lema\'{n}czyk \cite{Lemanczyk}, Hulse \cite{Hulse2}]\label{prop:sequenceQuasiProduct}
Let $(X,\mathscr{B},\mu,T)$ and $(X',\mathscr{B}', \mu',T')$ be two dynamical systems.
Then
\begin{equation}\label{eq:zza}
h_{A,\mu\times\mu'}(T\times T')\leq h_{A,\mu}(T)+h_{A,\mu'}(T').
\end{equation}
\end{proposition}

\brm
In \cite{Kushnirenko} it is claimed that one has equality in \eqref{eq:zza}. Counterexamples to that were given in \cite{Lemanczyk} (see also \cite{Hulse2}) where the author showed that the inequality is strict.
\erm

\brm
As shown by Lema\'{n}czyk \cite{Lemanczyk}, the formula $h_{A,\mu}(T^k)=kh_{A,\mu}(T)$, $k\geq 2$ is not true in general.  This  answered a question proposed by Saleski \cite{Saleski}.
\erm

Recall that the sequence entropy with the sequence of integers is by definition equal to the classical entropy. It is natural to ask whether there exists a class of more general sequences for which there is a relation between the sequence entropy and classical entropy. This was shown to be true for large density sequences and systems of positive (classical) entropy.  More precisely,  E.\ Krug and D.\ Newton have the following result to describe the relation between the sequence entropy and the classical entropy:

\begin{definition}[Krug-Newton \cite{KrugNewton}]\label{K(A)}
For a given  sequence of integers $A=\{t_1,\ldots\}$ and $k,n\in\mathbb{N}$, let $U_A(n,k)=\{t_i+j:1\leq i\leq n, -k\leq j\leq k\}$. Let $S_A(n,k)=\operatorname{Card} U_A(n,k)$ and define \footnote{Notice that $S_A(n,k)$ is an increasing positive function of $k$, that means $\overline{\lim}_{n\to\infty}\frac{1}{n}S_A(n,k)$ is an increasing non-negative  function of $k$. By allowing the value $\infty$, $K(A)$ is well defined.},
\begin{equation}
K(A):=\lim_{k\to\infty}\overline{\lim_{n\to\infty}}\frac{1}{n}S_A(n,k).
\end{equation}
\end{definition}

\begin{theorem}[Krug-Newton \cite{KrugNewton}]\label{thm:relationSequenceAndGeneralMetric}
Let $T$ be any automorphism of $(X, \mu)$ and $A$ any
sequence of integers. Then $h_{A,\mu}(T) = K(A)\cdot h(T)$,
where the right hand side is to be interpreted as $0$ if
$K(A)=0$ and $h(T)=\infty$ and is undefined if $K(A)=\infty$ and $h(T)=0$.
\end{theorem}

\brm
Notice that the above theorem does not provide any information for systems of zero entropy, for which one can hope to take sparser sequences to recover some information. For example, it was shown by Kushnirenko \cite{Kushnirenko} that the time one map $u_1$ of the horocycle flow satisfies: $1\leq h_{2^n,\mu}(u_1)\leq 6$ and this result was improved in \cite{KanigowskiVinhageWei} to $h_{2^n,\mu}(u_1)=3$. We provide more details in Section \ref{sec:sequenceEntropyQuasiU}.
\erm

\subsection{Spectral theory and measure-theoretic sequence entropy}\label{sec:metricSequenceEntropySpec}
In this section, we will discuss the relation between the spectrum of $T$ and the measure-theoretic sequence entropy of $T$.

We  start from the following result of Kushnirenko:
\begin{theorem}[Kushnirenko \cite{Kushnirenko}]\label{thm:totalVanishingSeqEnt}
An invertible transformation $T$ has discrete spectrum if and only if $h_{A,\mu}(T)=0$ for any sequence $A$.
\end{theorem}

\brm
This theorem should be compared with the corresponding result for the slow entropy in Section \ref{sec:CriterionTopMea}.
\erm




On the other hand we have the following result for systems with continuous spectrum:

\begin{theorem}[Pitskel \cite{Pitskel}]\label{thm:Pitskel}
Let $\xi$ be a $k$-element partition of the space $(X,\mu)$ and let $\chi_i$ be the characteristic function of the $i-$th element, $i=1,\ldots,k$. Let  $U_T$ be the Koopman operator and set $\psi_i:=\chi_i-\int\chi_id\mu$. Let $\psi_i$, $i=1,\ldots,k$ belong to an invariant subspace $L\subset L^2(X)$ and let the operator $U_T$ have a continuous spectrum on $L$. Then there exists a sequence $A$ such that
\begin{equation}
h_{A,\mu}(T,\xi)=H(\xi).
\end{equation}
\end{theorem}



Recall that for the classical entropy we have the generator theorem.
This turns out not to be correct for the sequence entropy. By strengthening Pitskel's result, P.\ Hulse obtained the following result:
\bt[Hulse \cite{Hulse}]
There exists a weakly mixing automorphism $T$ and an increasing sequence of natural numbers $A$ such that $h_{A,\mu}(T)=\infty$ and $h_{A,\mu}(T,\xi)=H(\xi)$ for all $\xi$ with finite entropy. Moreover, for any $\epsilon>0$, there exists a two-elements $A-$generator $\xi$ for $T$ such that $H(\xi)<\epsilon$.
\et

\subsection{Sequence entropy and spectral properties}
In this section we will describe further results on connections of the sequence entropy and spectral properties  (recall that first such result was established by Kushnirenko, Theorem \ref{thm:totalVanishingSeqEnt}).
The first two results describe how sequence entropy characterizes weak mixing:
\bt[Hulse \cite{Hulse}]
An invertible measure preserving system $(X,\mathscr{B},\mu,T)$ is weakly mixing if and only if there exists an increasing sequence of natural numbers $A$ such that $h_{A,\mu}(T,\xi)=H(\xi)$ for all $\xi$ with $H(\xi)<\infty$.
\et


The following result gives yet another characterization of weak mixing:

\bt[Zhang \cite{Zhang}]
An invertible measure preserving system $(X,\mathscr{B},\mu,T)$ is weakly mixing if and only if for any set $\Gamma\subset\mathbb{N}$ with positive density, there is a subset $\Gamma_1$ of $\Gamma$ such that for any finite partition $\xi$, $h_{\Gamma_1,\mu}(\xi)=H(\xi)$.
\et

The following results describe the connections between the sequence entropy, mild mixing, mixing and rigidity:
\bt[Zhang \cite{Zhang}]
An invertible measure preserving system $(X,\mathscr{B},\mu,T)$ is rigid if and only if there exists a subset $\Gamma$ of $\mathbb{N}$ such that if $\{F_i\}$ is any sequence of pairwise disjoint finite subsets of $\Gamma$ and $s_i=\sum_{a\in F_i}a$, then $h_{\{s_i\},\mu}(T)=0$.
\et

\bt[Zhang \cite{Zhang}]
An invertible measure preserving system $(X,\mathscr{B},\mu,T)$ is mildly mixing if and only if for any subset $\Gamma$ of $\mathbb{N}$, there is a sequence $\{F_n\}$ of pairwise disjoint finite subsets of $\Gamma$ such that for any finite partition $\xi$ of $X$ and $s_i=\sum_{a\in F_i}a$, $h_{\{s_i\},\mu}(T,\xi)=H(\xi)$.
\et

\bt[Saleski \cite{Saleski}, Hulse \cite{Hulse}]
An invertible measure preserving dynamical system $(X,\mathscr{B},\mu,T)$ is mixing if and only if for any infinite sequence $A$ there exists a subsequence $B\subset A$ such that $h_{B,\mu}(T,\xi)=H(\xi)$ for all $\xi$ with $H(\xi)<\infty$.
\et

\subsection{Topological sequence entropy}
In $1974$, seven years after  Kushnirenko's definition of measure-theoretic sequence entropy, Goodman \cite{Goodman} introduced topological sequence entropy and investigated its properties. Later, Huang, Shao and Ye \cite{HuangSahoYe} studied relations of topological weak mixing and topological sequence entropy which is a counterpart of the corresponding results in the measurable category from the previous section.We start this section with  the definition of the topological sequence entropy.


\begin{definition}[Topological sequence entropy, Goodman \cite{Goodman}]
Suppose that $X$ is a locally compact metric space and $T:X\to X$ is a continuous map. Let $\mathscr{A}(X)$ be the collection of all open covers of $X$. For a compact subset $K\subset X$ and a given open cover $\alpha$, let $N(\alpha,K)$ denote the minimal cardinality of a sub-cover of $\alpha$ that covers  $K$. For a given sequence of integers $A=\{t_1,\ldots,t_n,\ldots\}$, we define
\begin{equation}
\begin{aligned}
&h_A(T,\alpha,K):=\overline{\lim_{n\to\infty}}\frac{1}{n}H(\bigvee_{i=1}^nT^{-t_i}\alpha),\\
&h_A(T):=\sup_K\sup_{\alpha\in\mathscr{A}(X)}h_A(T,\alpha,K).
\end{aligned}
\end{equation}

\brm
If $X$ is a compact metric space we will omit the $K$ in the $h_A(T,\alpha,K)$ for simplicity.
\erm

The quantity $h_A(T)$ denotes the topological sequence entropy of $(X,T)$ with respect to the sequence $A$.
\end{definition}

If $T$ is a homeomorphism, then the above definition of topological sequence entropy can be extended to negative integers. Analogously to the  measure-theoretic sequence entropy, by taking $t_i=i-1$,   $h_A(T)$ coincides with the classical topological entropy.

Recall that the classical topological entropy can also be defined through maximal separated sets (or minimal spanning sets) with respect to the Bowen distance. There is also a parallel version of the topological sequence entropy's definition based on maximal separated sets (or minimal spanning sets).

\begin{definition}[Topological sequence entropy, Goodman \cite{Goodman}]
Suppose that $X$ is a locally compact metric space, $T:X\to X$ is a continuous map and $A=\{t_1,\ldots,t_n,\ldots\}$ is a given sequence of integers. For $\alpha\in\mathscr{A}(X)$, we say a set $E\subset X$ is $(A,n,\alpha)-$separated (with respect to $T$) if for any distinct points $x,y\in E$, there is an integer $j$, $1\leq j\leq n$ such that $T^{t_j}x\in U\in\alpha$ and $T^{t_j}y\notin U$. For a compact subset $K\subset X$ let $N(A,n,\alpha,K)$ denote the largest cardinality of a $(A,n,\alpha)-$separated set in $K$. Define
\begin{equation}
\begin{aligned}
&h_A(T,\alpha,K)=\overline{\lim_{n\to\infty}}\frac{1}{n}\log N(A,n,\alpha,K),\\
&h_A(T)=\sup_K\sup_{\alpha\in\mathscr{A}(X)}h_A(T,\alpha,K).
\end{aligned}
\end{equation}

The quantity $h_A(T)$ denotes  the topological sequence entropy of $(X,T)$ with respect to the sequence $A$.
\end{definition}
\begin{remark}
As we have already mentioned before, topological sequence entropy can also be defined by minimal spanning sets. We say a set $F$  $(A,n,\alpha)-$spans another set $K$ if for each point $x\in K$, there is a point $y\in F$ and there exists $U\in\alpha$ with $T^{t_j}x,T^{t_j}y\in U$ for some $j\leq n$. For a compact set $K\subset X$, let $S(A,n,\alpha,K)$ denote the smallest cardinality of a set which $(A,n,\alpha)-$spans $K$ and denote $S(A,\alpha,K)=\overline{\lim}_{n\to\infty}\frac{1}{n}\log S(A,n,\alpha,K)$. In this situation, define topological sequence entropy as
\begin{equation}
h_A(T,K):=\sup_{\alpha\in\mathscr{A}(X)}S(A,\alpha,K).
\end{equation}
\end{remark}

\brm
In fact, the  two above definitions are equivalent. Indeed, one can prove that
\begin{equation}
S(A,n,\alpha,K)\leq N(A,n,\alpha,K)\leq S(A,n,\gamma,K),
\end{equation}
where $\gamma$ is a refined open cover of $\alpha$.
\erm
\brm
Similar as for the measure-theoretic sequence entropy, Lema\'{n}czyk \cite{Lemanczyk} showed that topological sequence entropy also does not depend monotonically on $A$: it is possible that $h_A(T)<h_B(T)$  if $A$ is a subsequence of $B$.
\erm

Topological sequence entropy has also the following useful properties:

\begin{proposition}[Generator-type property, Goodman \cite{Goodman}]
 Let $\{\alpha_n\},n=1,\ldots$ be a sequence of elements of $\mathscr{A}(X)$ for which $diam(\alpha_n)\to0$. Then $h_A(T,\alpha_n)\to h_A(T)$ as $n\to\infty$.
\end{proposition}

\begin{proposition}[Cartesian square, Goodman \cite{Goodman}]
We have
\begin{equation}
h_A(T\times T)=2h_A(T).
\end{equation}
\end{proposition}

\begin{remark}
The above equality is not true for general products of two systems as shown by Lema\'{n}czyk \cite{Lemanczyk} (this answered a question of Goodman, \cite{Goodman}). One only has the following inequality:
\begin{equation}
h_A(T\times S)\leq h_A(T)+h_A(S).
\end{equation}
\end{remark}

\brm
As for iterates of $T$, Lema\'{n}czyk \cite{Lemanczyk} showed that it is not true that $h_A(T^k)=kh_A(T)$, $k\geq2$ which answered a question by Saleski \cite{Saleski} in the topological category.
\erm


The following two results describe the  behavior of the topological sequence entropy under  factors or extensions.
\begin{proposition}[Goodman \cite{Goodman}]
Let $\pi:(X,T)\to(Y,S)$ be a homomorphism and suppose each fiber $\pi^{-1}y$, $y\in Y$, has at most $n$ points. Let $A$ be any sequence. Then, $h_A(T)\leq h_A(S)+\log n$.
\end{proposition}

\begin{proposition}[Goodman \cite{Goodman}]
For any surjective dynamical system $(X,T)$ and sequence $A$, we have $h_A(T)=h_A(T^*)$, where $(X^*,T^*)$ is the natural extension of the $(X,T)$.
\end{proposition}

\subsection{Relations of topological sequence entropy and topological properties}\label{sec:TopologicalSeqMixing}

We first introduce some notation in the topological sequence entropy:
\bd
Let $(X,T)$ be a topological system.
\begin{enumerate}[(1)]
\item $(X,T)$ is called bounded if $h_A(T)<K$ for all sequences $A$, otherwise, we will say $(X,T)$ is unbounded.
\item $(X,T)$ is called null if the topological sequence entropy is zero for any sequence.
\item Let $S^1$ be the unit circle in the complex plane and $C(X)$ be the collection of all continuous maps from $X$ into $K$. Given a system $(X,T)$, we define a sequence of groups,
$$K=G_1\subset G_2\subset\ldots\subset C(X)$$
inductively as follows:
$$G_{n+1}=\{f\in C(X):f(Tx)=gf(x) \text{for all $x\in X$ and some $g\in G_n$}\}.$$
Denote by $n(T)$ the least number $n$ for which $G_{n+1}=G_n$.
\item An admissible cover $\mathcal{U}$ is a finite cover $\mathcal{U}=\{U_1,U_2,\ldots,U_n\}$ such that $\left(\bigcup_{j\neq i}U_j\right)^c$ has nonempty interior for each $i\in\{1,2,\ldots,n\}$;
\end{enumerate}
\ed

At first, a counterexample to the topological version of Kushnirenko's discrete spectrum theorem is described. Recall that by Theorem \ref{thm:totalVanishingSeqEnt}, vanishing of measure-theoretic sequence entropy for any sequence is equivalent to discrete spectrum. It is natural to conjecture a corresponding topological version, i.e.:
\begin{center}
A minimal flow is null if and only if it is equicontinuous.
\end{center}
In fact, this is not the case for topological sequence entropy. In \cite{Goodman}, Goodman gave an example which is a minimal null topological dynamical system which is not distal\footnote{Recall that equicontinuity implies distality.}.

\begin{example}[Goodman \cite{Goodman}]
Let $\alpha\in[0,1]$ be an irrational number and define $A_0=\{e^{2\pi i\theta}:0\leq\theta\leq\frac{1}{2}\}$, $A_1=\{e^{2\pi i\theta}:\frac{1}{2}\leq\theta\leq1\}$ and $X=\{x\in\{0,1\}^{\mathbb{Z}}:\cap_{r=-\infty}^{\infty}e^{2\pi ir\alpha}A_{x_r}\neq\emptyset\}$. We also denote $T$ as the shift given by $(Tx)_n=x_{n+1}$ for $x\in X$. Then $(X,T)$ is a minimal null topological system which is not distal.
\end{example}

Moreover, the nullness property is quite stable under common operations and also related to uniquely ergodicity.

\bt[Huang-Li-Shao-Ye \cite{HLSY}]
The property of nullness of a topological dynamical system is stable under factor maps, arbitrary products and inverse limits. Moreover, if a null system is minimal, then it is uniquely ergodic.
\et
We also have the following characterization of topological weak mixing:
\bt[Goodman \cite{Goodman}]
Let $(X,T)$ be an invertible topological dynamical system:
\begin{enumerate}[(1)]
\item If $(X,T)$ is non-trivial and topologically weakly mixing, then $(X,T)$ is unbounded;
\item If $(X,T)$ is minimal, then it is topologically weakly mixing if and only if it has no non-trivial bounded factors;
\item Assume that $(X,T)$ is a non-trivial bounded system. If $(X,T)$ is minimal, then $n(T)\geq2$. If $(X,T^r)$ is minimal for all integers $r\neq 0$, then $n(T)\leq2$.
\end{enumerate}
\et

With the help of admissible covers we can describe topological mixing, topological weak mixing and topological mild mixing as follows:
\bt[Huang, Shao and Ye \cite{HuangSahoYe}] Let $(X,T)$ be a topological dynamical system, then we have
\begin{enumerate}[(1)]
  \item $(X,T)$ is topologically weakly mixing iff for any admissible open cover $\mathscr{U}$ of $X$ there is some $A\subset\mathbb{Z}_+$ such that $h_A(T,\mathscr{U})=\log N(\mathscr{U})$, where $N(\mathscr{U})$ is the cardinality of the $\mathscr{U}$;
  \item (X,T) is topologically mildly mixing iff for any admissible open cover $\mathscr{U}$ if $X$ and any IP set $F$ there is some $A\subset F$ such that $h_A(T,\mathscr{U})=\log N(\mathscr{U})$, where IP set $F$ is defined as $$F=\{\sum_{i\in\alpha}b_i:\text{$\alpha$ is a finite non-empty subset of $\mathbb{N}$}\},$$ for a sequence of natural numbers $b_i$;
  \item Assume that $(X,T)$ is topologically mixing. Then for any admissible open cover $\mathscr{U}$ of $X$ and any infinite $F\subset\mathbb{Z}_+$ there is some infinite $A\subset F$ such that $h_A(T,\mathscr{U})=\log N(\mathscr{U})$.
\end{enumerate}
\et

\subsection{Variational principal for sequence entropy}
Recall that  the connection between classical topological and measure theoretic entropy is given by the variational principle. It is therefore interesting to ask for an analogous property for the sequence entropy. We will discuss this in this section.

Recall that $\mathscr{M}(X,T)$ denotes the collection of all regular probability measure on Borel subsets of $X$ that are invariant under the $T$.
\bd[Finite covering dimension, Goodman \cite{Goodman}]
A space $X$ has finite covering dimension at most $n$ if any open cover of $X$ has an open refinement of order at most $n+1$. The order of a collection $\alpha$ of sets is defined to be the maximum number of sets in $\alpha$ with non-empty intersection.
\ed
\begin{theorem}[Goodman \cite{Goodman}]
Let $X$ be a compact Hausdorff space and $T$ be a continuous map from $X$ to $X$. Suppose $X$ has finite covering dimension. Then for any sequence $A$ and any measure $\mu\in \mathscr{M}(X,T)$,
\begin{equation}
h_{A,\mu}(T)\leq h_A(T).
\end{equation}
\end{theorem}

The opposite inequality also holds under some additional assumptions (see Theorem \ref{thm:relationSequenceAndGeneralMetric} by Krug and Newton):

\begin{theorem}[Goodman \cite{Goodman}]
Let $X$ be a compact Hausdorff space and $T$ be a continuous map from $X$ to $X$. Suppose $X$ has finite covering dimension. Then for any sequence $A$ suppose that either $K(A)<\infty$ or $h_{A}(T)>0$, then,
\begin{equation}
h_A(T)=\sup_{\mu\in\mathscr{M}(X,T)}h_{A,\mu}(T).
\end{equation}
\end{theorem}

\begin{remark}
In \cite{Goodman}, a counterexample is provided to show that without the condition in Theorem \ref{thm:relationSequenceAndGeneralMetric}, we may have\footnote{More precisely, the example is a Chacon system which has a special spacers sequence.} $$h_A(T)>\sup_{\mu\in\mathscr{M}(X,T)}h_{A,\mu}(T).$$
\end{remark}


\subsection{Sequence entropy of quasi-unipotent systems}\label{sec:sequenceEntropyQuasiU}
In this subsection, one special example related to sequence entropy is introduced, i.e.\ the sequence entropy of quasi-unipotent flows. This example is a generalization  of the original Kushnirenko's computation for  horocycle flows.

\subsubsection{Introduction to quasi-unipotent flows}\label{sec:introQuasiUni}
Let $G$ be a connected Lie group and $\mathfrak{g}$ its Lie algebra. Suppose $\Gamma$ is a discrete subgroup of $G$ with co-finite volume and $\mu$ is the Haar measure on $G\slash \Gamma$ which is induced from Riemannian volume on $G$. We define a quasi-unipotent flow as follows:
\bd[Quasi-unipotent flow]
An element $U \in \mathfrak{g}$ is called {\it quasi-unipotent} if it can be written as $U = U' + Q$, where $\mathrm{ad}_{U'}^N = 0$ for some $N$, $Q$ is $\mathrm{ad}$-compact, and $[Q,U'] =0$. A quasi-unipotent element $U \in \mathrm{Lie}(G)$ induces a quasi-unipotent flow on a homogeneous space $G\slash \Gamma$ by

$$\varphi_t(g\Gamma) = \exp(tU)g \Gamma.$$
\ed

In order to describe the complexity of quasi-unipotent flows, the following structure is introduced as a basis system for homogeneous dynamics.

\begin{definition}
	\label{def:chain-basis}	
	A {\em chain} in $\mathfrak{g}$ with respect to a quasi-unipotent element $U$ of depth $m$ is a linearly independent set $\{ X_j : 0 \leq j \leq m\}$  such that $X_0$ is in the centralizer of $U$ and
	
	\[ \mathrm{ad}_U(X_j) = X_{j-1} \mbox{ for all }1 \leq j \leq m.\]

 A {\em double chain} in $\mathfrak{g}$ with respect to $U$ of depth $m$ is a linearly independent set \\ $\{X_{j,i} : 0 \leq j \leq m, i = 0,1\}$ and a number $\alpha$ such that $\mathrm{ad}_Q(X_{j,0}) = -\alpha X_{j,1}$, $\mathrm{ad}_Q(X_{j,1}) = \alpha X_{j,0}$ for all $0 \leq j \leq m$, $X_{0,i}$ is in the centralizer of $U'$ for $i = 0,1$ and

\[ \mathrm{ad}_{U'}(X_{j,i}) = X_{j-1,i} \mbox{ for all }1 \leq j \leq m, i = 0,1.\]

	A {\em chain basis} of $\mathfrak{g}$ with respect to $U$ is a basis of chains and double chains. The sequence of depths $(m_1,\dots,m_n)$ of chains and double chains, with each double chain listed twice, is called the {\em chain structure} of $U$.

\end{definition}

It follows from the Jordan block decomposition that all quasi-unipotent elements have a chain basis and thus for a quasi-unipotent flow $\varphi_t(g\Gamma)= \exp(tU)g\Gamma$ with chain structure $(m_1,m_2,\ldots,m_n)$, we define
$$GR(U)=\sum_{i=1}^n\frac{m_i(m_i+1)}{2}.$$

\subsubsection{Sequence entropy of quasi-unipotent flows}
The sequence entropy of quasi-unipotent flows is described in the following theorem:

\begin{theorem}[Kanigowski-Vinhage-Wei \cite{KanigowskiVinhageWei}]\label{thm:mainco}
For a quasi-unipotent flow $\varphi_t$ on homogeneous space $G\slash \Gamma$ with chain structure $(m_1,m_2,\ldots,m_n)$, Haar measure $\mu$ and the sequence $A=\{C\lambda^n\}$ ($C>0$), we have
\begin{equation}
h_A(\varphi_t)=h_{A,\mu}(\varphi_t)=GR(U) \log \lambda.
\end{equation}
\end{theorem}

\brm
This result extends Kushnirenko's original result in \cite{Kushnirenko} by calculating the precise value of the sequence entropy of quasi-unipotent flows. For example, the sequence entropy\footnote{By taking logarithm function with base $2$, which is the same as in Kushnirenko's original paper \cite{Kushnirenko}.} of horocycle flows  with respect to the sequence $\{2^n\}$ is equal to $3$.
\erm
\brm
It is also worth to point out that the topological sequence entropy equals the measure-theoretic sequence entropy for the Haar measure, i.e. the variational principle holds for quasi-unipotent flows, for more details see \cite{KanigowskiVinhageWei}.
\erm

\section{Slow entropy}\label{sec:slowMetricEntropy}
In this section, a basic theory of slow entropy is introduced, applications of slow entropy to smooth realization and spectrum are discussed. Moreover, several special features of slow entropy are also discussed such as failure of variational principle and difference between upper and lower definitions of slow entropy. Many examples of slow entropy are also studied: quasi-unipotent flows, rank one systems, surface flows and AbC constructions.\subsection{Measure-theoretic slow entropy}\label{sec:SlowMeaureTHeoreticEntropy}
In this subsection following \cite{KatokThou}, the definition of the slow entropy is introduced in the amenable group setting.

\bd
Suppose that $\Gamma$ is an amenable discrete group and $F$ is a subset of $\Gamma$. Then
$$\Omega_{N,F}=\{w=(w_r)_{r\in F};w_r\in\{1,\ldots,N\}\}$$
and define a natural projection $\pi_{F,F'}:\Omega_{N,F}\to\Omega_{N,F'}$ for $F'\subset F$.
\ed

Suppose $T:(X,\mu)\times\Gamma\to(X,\mu)$ is an action of the group $\Gamma$ by measure preserving transformations of a Lebesgue space. Let $\xi=\{c_1,\ldots,c_N\}$ be a finite measurable partition. We define the coding map as $\phi_{T,\xi}:X\to\Omega_{N,\Gamma}$ if $(\phi_{T,\xi})_r(x)=w_r(x)$ where $T(r)x\in c_{w_r(x)}$. Denote the partial coding  $\phi^F_{T,\xi}$ for $F\subset\Gamma$ as $\phi^F_{T,\xi}=\pi_{\Gamma,F}\circ\phi_{T,\xi}$ and call $\phi^F_{T,\xi}$ the $F-$ {\em name of $x$ with respect to $\xi$}. The measure induced by a partial coding is $(\phi^F_{T,\xi})_*\mu$. Then, the Hamming metric on $\Omega_{N,F}$ is defined as follows.

\begin{definition}[Hamming metric]\label{def:HammingMetric}
For any finite set $F\subset\Gamma$, the Hamming metric $d_F^H(w,w')$ is defined by
 \begin{equation}
d^H_F(w,w')=\frac{1}{\mathrm{card} \,F}\sum_{r\in F}(1-\delta_{\omega_r\omega'_r}).
\end{equation}
\end{definition}

Now the cover of symbolic space $\Omega_{N,F}$ can be introduced based on the Hamming metric. Denote by $S_{\xi}^H(T,F,\epsilon,\delta)$ the minimal number of $d_F^H-\epsilon-$balls in $\Omega_{N,F}$, whose union has measure at least $1-\delta$ with respect to $(\phi^F_{T,\xi})_*\mu$. Finally we assume that $\{F_n\}_{n=1}^{\infty}$ is a F{\o}lner sequence of finite subsets of $\Gamma$. Then the slow entropy can be defined as:

\begin{definition}[Slow entropy, Katok-Thouvenot \cite{KatokThou}]\label{def:slowEntropy}
Let ${\bf a}=\{a_n(t)\}_{n\in\mathbb{N},t>0}$ be a family of positive sequences increasing to infinity and monotone in $t$. Then define the upper measure-theoretic slow entropy of $T$ with respect to $\xi$ by
\begin{equation}\label{eq:slowEntroDef}
\overline{ent}_{\bf a}^{\mu}(T,\xi)=\lim_{\delta\to0}\lim_{\epsilon\to0} A(\epsilon,\delta,\xi),
\end{equation}
where
$
A(\epsilon,\delta,\xi)=\left\{
\begin{array}{ll}
\sup B(\epsilon,\delta,\xi), & \hbox{if $B(\epsilon,\delta,\xi)\neq\emptyset$;} \\
0, & \hbox{if $B(\epsilon,\delta,\xi)=\emptyset$,}
\end{array}
\right.
$
for $$B(\epsilon,\delta,\xi)=\{t>0:\limsup_{n\to\infty}\frac{S_{\xi}^H(T,F_n,\epsilon,\delta)}{a_n(t)}>0\}.$$

The upper measure-theoretic slow entropy of $T$ is defined by
\begin{equation}
\overline{ent}_{\bf a}^{\mu}(T)=\sup_{\xi}\overline{ent}_{\bf a}^{\mu}(T,\xi).
\end{equation}
\end{definition}

\brm
Similarly, the lower $a-$entropy can be defined by changing the $\limsup$ to $\liminf$ in the above and use $\underline{ent}_{\bf a}^{\mu}(f)$ to represent the lower $a-$entropy of $T$.
\erm
\brm
By considering the continuous version of Hamming metric, Definition \ref{def:slowEntropy} can be modified to  flows.
\erm

If $T$ is ergodic and the scaling function is exponential: $a_n(t)=\exp(tn)$, then the slow entropy invariant coincides with the classical measure-theoretic entropy. The following theorem explains this phenomenon and also should be understood as a general version of Theorem \ref{MetricEntropyBall} in \cite{KatokPubl} and Proposition $2$ in \cite{Ferenczi}.
\bt[Katok-Thouvenot \cite{KatokThou}]\label{SlowEntropyGeneralEntropy}
Let $\{F_n\}$ be a F\o lner sequence in $\Gamma$. If $T$ is ergodic, then
$$\lim_{n\to\infty}\lim_{\epsilon\to0}\lim_{\delta\to0}\frac{\log S_{\xi}^H(T,F_n,\epsilon,\delta)}{|F_n|}=h_{\mu}(T).$$
\et


Recall that in the classical measure-theoretic entropy setting, a very powerful theorem in calculation of entropy is the generator theorem (Theorem \ref{GeneratorTheorem}). In fact, this is also the case for the slow entropy:

\bt[Katok-Thouvenot \cite{KatokThou}]\label{thm:GeneratorTheoremSlow}
If $\xi$ is a generator then $\overline{ent}_a^{\mu}(T)=\overline{ent}_a^{\mu}(T,\xi)$.
More generally, if $\xi_m$ is a sufficient family of partitions, then
$$\overline{ent}^{\mu}_a(T)=\sup_m\overline{ent}^{\mu}_a(T,\xi_m).$$
\et

The upper and lower slow entropy have different behavior under products:
\bp[Ferenczi \cite{Ferenczi}]\label{prop:slowEntropyProduct} If the scaling function satisfies $a_n(t+s)=a_n(t)a_n(s)$, we have:
\beq
\begin{aligned}
&\underline{ent}_a^{\mu}(T_1\times T_2)\geq\underline{ent}_a^{\mu}(T_1)+\underline{ent}_a^{\mu}(T_2),\\
&\overline{ent}_a^{\mu}(T_1\times T_2)\leq\overline{ent}_a^{\mu}(T_1)+\overline{ent}_a^{\mu}(T_2).
\end{aligned}
\ee
\epo




At the end of this section, several natural examples are provided to calculate slow entropy:

\begin{example}
Consider the Morse system $\sigma_M$ generated by the following substitution:
\begin{equation}
\begin{aligned}
a\to ab,\\
b\to ba.
\end{aligned}
\end{equation}
By calculations in \cite{Mosse}, it follows that the upper growth of the system has rate $\frac{10n}{3}$ and lower growth rate is $3n$ which means they are of the same (linear) order.

Thus we have $$\overline{ent}_a^{\mu}(\sigma_M)=\underline{ent}_a^{\mu}(\sigma_M)=1,$$ where $a_n(t)=n^t$.
\end{example}

\begin{example}
Consider the Rudin-Shapiro system $\sigma_{RS}$ generated by the following substitution:
\begin{equation}
\begin{aligned}
a\to ab,\\
b\to ac,\\
c\to db,\\
d\to dc.
\end{aligned}
\end{equation}
Then this system upper growth and lower growth coincide and the asumptotic value is $8n$. Thus we have,
$$\overline{ent}_a^{\mu}(\sigma_{RS})=\underline{ent}_a^{\mu}(\sigma_{RS})=1,$$ where $a_n(t)=n^t$.
\end{example}
\subsection{Topological slow entropy and Goodwyn's theorem}\label{sec:SlowTopologicalEntropy}
 We will follow the ideas from \cite{KatokThou} and \cite{Galatolo} to define the topological slow entropy. Suppose $K \subset X$ is a compact subset of a locally compact metric space $(X,d)$ and $T:(X,d)\times\Gamma\to(X,d)$ is an action of the group $\Gamma$ by homeomorphisms. Let $F$ be a finite subset of $\Gamma$ and we define Bowen metric as $d_F^T=\max_{\gamma\in F}d\circ T(\gamma)$. Moreover, we assume that $\{F_n\}_{n=1}^{\infty}$ is a F{\o}lner sequence of finite subsets of $\Gamma$. We will modify the definition in Section \ref{sec:topologicalentropy} a little to cover the case where $X$ is non-compact.

\bd[Minimal separated set and maximal spanning set]
Let $N_{d_F^T,K}(\epsilon)$ be the minimal number of $d_F^T-\epsilon-$balls required to cover $K$ (since $K$ is compact, this is finite);
Let $S_{d_F^T,K}(\epsilon)$ be the maximal number of $d_F^T-\epsilon-$balls  with centers in $K$ which can be placed disjointly in $X$.
\ed

Then the topological slow entropy can be defined as follows:
\begin{definition}[Topological slow entropy]
Let $\{a_n(t)\}_{n\in\mathbb{N},t>0}$ be a family of positive sequences increasing to infinity and monotone in $t$. The topological slow entropy of $T$ with respect to $\{a_n(t)\}_{n\in\mathbb{N},t\geq0}$ is
\begin{equation}
\overline{ent}^{\operatorname{top}}_{a}(T) =\sup_K \lim_{\epsilon \to 0} N(\epsilon,K),
\end{equation}
where $N(\epsilon,K)=\left\{\begin{array}{ll}
 \sup N_1(\epsilon,K), & \hbox{if $N_1(\epsilon,K)\neq\emptyset$;} \\
0, & \hbox{if $N_1(\epsilon,K)=\emptyset$,}
\end{array}
\right.$
for $$N_1(\epsilon,K)=\{t>0: \limsup_{n \to \infty} \dfrac{N_{d_{F_n}^T,K}(\epsilon)}{a_n(t)} > 0\}.$$

Equivalently, $\overline{ent}^{\operatorname{top}}_{a}(T)$ can also be defined as:
\begin{equation}
\overline{ent}^{\operatorname{top}}_{a}(T)=\sup_K \lim_{\epsilon \to 0} S(\epsilon,K),
\end{equation}
where
$S(\epsilon,K)=\left\{
\begin{array}{ll}
 \sup S_1(\epsilon,K), & \hbox{if $S_1(\epsilon,K)\neq\emptyset$;} \\
0, & \hbox{if $S_1(\epsilon,K)=\emptyset$,}
\end{array}
\right.$
for $$S_1(\epsilon,K)=\{ t>0: \limsup_{n \to \infty} \dfrac{S_{d_{F_n}^T,K}(\epsilon)}{a_n(t)} > 0\}.$$
\end{definition}
\brm
If we pick $a_n(t) = e^{nt}$, then topological slow entropy defined above coincides with the classical topological entropy.
\erm

\subsubsection{Slow entropy Goodwyn's theorem}
It is natural to ask for the variational principle for measure-theoretic and topological slow entropy. In fact, the answer is negative and counterexamples are provided in Section \ref{sec:failureVariational}. On the other hand, we still have one sided inequality in the slow entropy setting, i.e.\ slow entropy Goodwyn's theorem. Indeed, this is a corollary of Proposition $2$ in \cite{KatokThou}. We introduce the following definition to formulate the slow entropy Goodwyn's theorem.

\bd
A metric space $X$ is well-partitionable if for any Borel probability measure $\mu$, compact set $K$ and $\epsilon > 0$, there exists a finite partition $\mathscr{P}$ of $K$ whose atoms have diameter less than $\epsilon$ and such that $\mu\left(\bigcup_{\xi \in \mathscr{P}} \partial_\epsilon\xi\right) < \epsilon$, where $\partial_{\epsilon}\xi = \{ y \in X : B(y,\epsilon) \cap \xi \neq \emptyset \text{ but } B(y,\epsilon) \nsubseteq \xi\}.$
\ed

Notice that any smooth manifold is well-partitionable\footnote{In \cite{KatokThou} the authors consider a compact space $X$ but their proof of Proposition $2$ generalizes easily to the case where $X$ is  well-partitionable.}.

\begin{theorem}[Slow entropy Goodwyn's Theorem]\label{thm:Goodwyn}
Suppose $X$ is well-partitionable and $T:(X,d)\times\Gamma\to(X,d)$ is an action of the group $\Gamma$ by homeomorphisms preserving a non-atomic Borel probability measure $\mu$. Let  $\{a_n(t)\}_{n\in\mathbb{N},t>0}$ be a family of sequences increasing to infinity and monotone in $t$, then
$$ \overline{ent}^{\mu}_{a}(T)\le \overline{ent}^{\operatorname{top}}_{a}(T).$$
\end{theorem}

\subsection{Generalized Kushnirenko's inequality}
As we have already mentioned at the beginning of Section \ref{sec:KushnirenkoSection}, slow entropy is a useful tool in the smooth realization problem. In fact, this is also one of the main results of \cite{KatokThou}. In \cite{KatokThou}, the authors not only gave a generalized version of Kushnirenko's inequality in the slow entropy setting, but also constructed an example  to distinguish two versions of Kushnirenko's inequality, i.e. a $\mathbb{Z}^2-$action which cannot be smoothly realized but satisfies  Kushnirenko's inequality (Theorem \ref{thm:KushnirenkoInequality}) and fails the generalized Kushnirenko inequality (Theorem \ref{thm:NewKushnirenkoInequality}).

The generalized Kushnirenko inequality will be derived from the following proposition.

\bp[Katok-Thouvenot \cite{KatokThou}]\label{prop:NewKushnirenkoInequality}
Assume that the group $\Gamma$ is generated by a finite set $\Gamma_0$ and that the F{\o}lner sets $F_n$ consists of elements whose word-length norm with respect to $\Gamma_0$ does not exceed $a_n$. Suppose that the compact metric space $X$ has finite box dimension $D$ with respect to the metric $d$ and $\Gamma$ acts on $X$ by \textbf{bi-Lipschitz homeomorphisms}. Denote by $L$  a common Lipschitz constant for all elements $T(r)$, $r\in\Gamma_0\cup\Gamma_0^{-1}$. Then for any partition $\xi$ with $\mu(\partial\xi)=0$, and some constant $c(\epsilon)$ only depending on $\epsilon$ and $X$, we have
$$S_{\xi}^H(T,F_n,\epsilon,\epsilon)\leq c(\epsilon)\exp(Da_n\log L).$$
\epo

Combining the above proposition with Theorem \ref{thm:GeneratorTheoremSlow}, the generalized Kushnirenko inequality can be formulated as follows:

\bt[Generalized Kushnirenko's inequality, Katok-Thouvenot \cite{KatokThou}]\label{thm:NewKushnirenkoInequality}
Under the same assumptions as in Proposition \ref{prop:NewKushnirenkoInequality}, a measure-preserving action of $\Gamma$ by Lipschitz homeomorphisms in a finite box dimensional space can be smoothly realized if its slow entropy is no more than exponential measured by the diameter of the F\o lner set in the word-length metric.
\et

Notice that zero entropy of the action corresponds to the subexponential growth measured by the number of elements in the F\o lner set. Thus, for any group which is not a finite extension of $\mathbb{Z}$ there is a gap between the classical entropy and slow entropy. In fact, this indicates that zero entropy system with superexponential growth in terms of diameters of F\o lner sets can be constructed and such systems do not have smooth or Lipschitz realization.

\subsubsection{Failure of Kushnirenko's inequality}
Let $C_n=[0,n-1]\times[0,n-1]\subset\mathbb{Z}^2$ and $T_1$, $T_2$ are two commuting measure-preserving transformations generating the $\mathbb{Z}^2-$action $T$ on the space $(X,\mu)$. Then by \cite{KatokThou}, the existence of the following special $\mathbb{Z}^2-$action is guaranteed:

\bp[Katok-Thouvenot \cite{KatokThou}]\label{prop:failKushnirenko}
Let $\epsilon_n$, $n=1,2,\ldots$ be a sequence of positive numbers decreasing to $0$. There exists a measure-preserving ergodic $\mathbb{Z}^2-$action $T$ on the unit interval with Lebesgue measure generated by transformations $T_1,T_2$ and a generating partition $\xi$ such that

(1) for some $\epsilon>0$, $\delta>0$,
$$\limsup_{n\to\infty}\frac{\log(S_{\xi}^H(T,C_n,\epsilon,\delta))}{n^2\epsilon_n}\geq1;$$

(2) the $\mathbb{Z}^2$-entropy $h_{\mu}(T)=0$;

(3) every element of the $\mathbb{R}^2$ suspension action has zero entropy.
\epo

As a corollary of the above proposition we get the following:

\bc[Katok-Thouvenot \cite{KatokThou}]
If $\lim_{n\to\infty}n\epsilon_n=\infty$, then the action constructed in Proposition \ref{prop:failKushnirenko} is not isomorphic to any action by Lipschitz homeomorphisms on a compact metric space of finite box dimension preserving a Borel probability measure. In particular, it is not isomorphic to an action by diffeomorphisms of a compact differentiable manifold.
\ec

The main method used in the construction of the above action is the cutting and stacking method by alternating periodic and independent concatenations. Interested readers may go to \cite{KatokThou} for details.

Slow entropy can also be used for smooth realization problems in the setting of infinite-measure preserving systems. By generalizing slow entropy to infinite-measure preserving systems, M. Hochman \cite{Hochman} produced examples of infinite-measure preserving $\mathbb{Z}^2-$actions that are not isomorphic to actions by diffeomorphisms on a compact manifold preserving an infinite Borel measure.

\subsection{Criterion of total vanishing of slow entropy}\label{sec:CriterionTopMea}

Recall that for a dynamical system $(X,\mathscr{B},\mu,T)$, if its sequence entropy is zero for all subsequences of $\{1,2,\ldots,\}$, then it has discrete spectrum. As slow entropy is also an invariant which measures the complexity of a given system, it is natural to expect that the total vanishing of slow entropy systems may have some analogous spectral consequences. Indeed, by \cite{Ferenczi} (measure-theoretic slow entropy) and \cite{KanigowskiVinhageWei} (topological slow entropy), total vanishing of slow entropy at all scales will restrict the system to being isomorphic to a translation in respectively measurable and topological category. This is the content of the two results below:

\begin{theorem}[Ferenczi \cite{Ferenczi}]\label{thm:measurable-zero-all-scales}
$T$ is measure-theoretically isomorphic to a translation on a compact group if and only if
\begin{equation}
\overline{ent}_a^{\mu}(T)=0\text{ or }
\underline{ent}_a^{\mu}(T)=0
\end{equation}
with respect to every family of scales $a_n(t)$.
\end{theorem}

\begin{theorem}[Kanigowski-Vinhage-Wei \cite{KanigowskiVinhageWei}]
\label{thm:topological-zero-all-scales}
A minimal homeomorphism $T : X \to X$ of a compact metric space
is topologically conjugated to a translation on a compact abelian group if and only if $\overline{ent}^{\operatorname{top}}_{a}(f)=0$ for every family of scales $a_n(t)$.
\end{theorem}

\subsection{Failure of variational principle for slow entropy}\label{sec:failureVariational}
Variational principle is one of the most important result in classical entropy theory.  Unfortunately, it does not hold for slow entropy type invariants. Although we still have Goodwyn's theorem for slow entropy, there are examples for which the topological slow entropy may not be equal to the supremum of measure-theoretic slow entropies, even for uniquely ergodic systems.

A general idea for the construction of counterexamples is based on Theorem \ref{thm:measurable-zero-all-scales} and Theorem \ref{thm:topological-zero-all-scales}: One finds a minimal topological system $f : X \to X$ all of whose ergodic measures yield discrete spectrum systems (Kronecker systems) but which is not topologically conjugate to a translation on a compact abelian group. In fact, this is possible by the {\it approximation-by-conjugation method} first used by Anosov and Katok in \cite{AnosovKatok}. More precisely, we have:
\begin{proposition}
Let $M$ be a manifold with a free circle action. Then there exists a uniquely ergodic, volume preserving, $C^{\infty}$ diffeomorphism $f : M \to M$ which is measurably conjugate to a translation on a torus $\mathbb{T}^d$, $d \geq 1$.
\end{proposition}

If $f : X \to X$ is a homeomorphism of a compact metric space, let $\mathfrak{M}(f)$ denote the space of invariant measures. We may also use the example due to Furstenberg:

\begin{proposition}[Furstenberg \cite{Fur}]
There exists a minimal $C^{\infty}$ diffeomorphism of $\mathbb{T}^2$, with a 1-parameter family of measures all of which yield discrete spectrum systems and isomorphic to $R_\alpha$ for some fixed $\alpha \in S^1$.
\end{proposition}

As a result we establish the failure of variational principle for slow entropy:

\begin{corollary}[Kanigowski-Vinhage-Wei \cite{KanigowskiVinhageWei}]
Let $M$ be a manifold with a free circle action. Then there exists a $C^{\infty}$ diffeomorphism $f : M \to M$ and family of scales $a_{\chi}$ such that
$$ \sup_{\mu \in \mathfrak{M}(f)} \overline{ent}_{a_\chi}^{\mu}(f) < \overline{ent}^{\operatorname{top}}_{a_\chi}(f). $$
\end{corollary}

\subsection{Different upper and lower quantities of slow entropy}\label{sec:DifferentUpLow}

Notice in the definition of slow entropy we may take $\limsup$ or $\liminf$ in \eqref{eq:slowEntroDef} to get respectively $\overline{\operatorname{ent}}_a^{\mu}$ or $\underline{\operatorname{ent}}_a^{\mu}$, respectively. Recall that for classical entropy the two above quantities coincide. As follows by a result of Cyr and Kra \cite{CyrKra}, this is not the case for slow entropy:
\bt[Cyr-Kra \cite{CyrKra}]\label{thm:CyrKraDifferenceUpLow}
Assume $(a_n)_{n\in\mathbb{N}}$ and $(b_n)_{n\in\mathbb{N}}$ are two non-decreasing sequences of positive integers such that $\lim_{n\to\infty}a_n=\infty$, $\lim_{n\to\infty}\frac{1}{n}\log b_n=0$ and $a_n\leq b_n$ for all $n\in\mathbb{N}$. There exists a minimal subshift $(X_{\infty},\sigma)$ of topological entropy zero and a $\sigma$-invariant, ergodic measure $\mu$ supported on $X_{\infty}$ such that
\beq\label{eq:CyrKraResults}
\begin{aligned}
&\sup_{\xi}\lim_{\delta\to0}\lim_{\epsilon\to0}\liminf_{n\to\infty}\frac{S_{\xi}^H(\sigma,n,\epsilon,\delta)}{a_n}\leq1,\\
&\sup_{\xi}\lim_{\delta\to0}\lim_{\epsilon\to0}\limsup_{n\to\infty}\frac{S_{\xi}^H(\sigma,n,\epsilon,\delta)}{b_n}\geq1.
\end{aligned}
\ee
\et


Let
\beq
s_n(t):=\left\{
  \begin{array}{ll}
   b_n(2t-1), & \hbox{$1\leq t$;}\\
    a_n(2-2t)+b_n(2t-1), & \hbox{$\frac{1}{2}\leq t\leq1$;} \\
    2ta_n+(1-2t)\frac{a_0}{2}, & \hbox{$0\leq t\leq\frac{1}{2}$.}
  \end{array}
\right.
\ee
Due to $a_n\leq b_n$, $\lim_{n\to\infty}a_n=\infty$ and $\lim_{n\to\infty}\frac{1}{n}\log(b_n)=0$, it is clear that $\{s_n(t)\}_{n\in\mathbb{N},t>0}$ is a family of positive sequences increasing to infinity for $t>0$, monotone in $t$ for $t\in[0,\infty)$ and $s_n(0)$ is a positive constant sequence. Moreover, $s_n(t)$ satisfies $s_n(\frac{1}{2})=a_n$ and $s_n(1)=b_n$.

By using the language of slow entropy, Theorem \ref{thm:CyrKraDifferenceUpLow} can be reformulated as follows:
\bt
Let $(s_n)_{n\in\mathbb{N}}$ be defined as above. There exists a minimal subshift $(X_{\infty},\sigma)$ of topological entropy zero and an ergodic measure $\mu$ supported on $X_{\infty}$ such that
\beq
\underline{ent}_s^{\mu}(\sigma)\leq\frac{1}{2}<1\leq\overline{ent}_s^{\mu}(\sigma).
\ee
\et







\subsection{Slow entropy of quasi-unipotent flows}\label{sec:QuasiUnipotent}
Recall that quasi-unipotent flows have zero entropy. Therefore, it is a natural idea to use slow entropy to study quasi-unipotent flows. In this section, we mainly state results from \cite{KanigowskiVinhageWei}.
Define $G$ as a connected Lie group and $\mathfrak{g}$ as its Lie algebra. Suppose $\Gamma$ is a discrete subgroup of $G$ with co-finite volume and $\mu$ is the Haar measure on $G/\Gamma$which is induced from Riemannian volume on $G$.

\begin{theorem}[Kanigowski-Vinhage-Wei \cite{KanigowskiVinhageWei}]\label{thm:slowentropyinvariants}
\label{thm:main}
Let $\varphi_t$ be an ergodic quasi-unipotent flow\footnote{See Section \ref{sec:introQuasiUni} for quasi-unipotent flow's definition.} on a homogeneous space $G/\Gamma$ with chain structure $(m_1,\dots,m_n)$. Then the topological and measure-theoretic slow entropy of $\varphi_t$ with respect to $a_n(\chi)=n^{\chi}$ coincide,
$$\overline{ent}^{\mu}_{a}(\varphi_t)=\underline{ent}^{\mu}_{a}(\varphi_t)=\overline{ent}^{\operatorname{top}}_{a}(\varphi_t)=\underline{ent}^{\operatorname{top}}_{a}(\varphi_t)=GR(U),$$
where $GR(U)$ is defined Section \ref{sec:introQuasiUni}.
\end{theorem}
\brm
As a corollary to Theorem \ref{thm:main}, one obtains a version of the variational principle for quasi-unipotent flows and slow entropy. Namely, the topological slow entropy of a quasi-unipotent flows is the supremum of the measure-theoretic slow entropies taken over all invariant measures. Recall that by Section \ref{sec:failureVariational} this property is fairly special.
\erm

\subsubsection{Examples of quasi-unipotent flows}
We give an example of a quasi-unipotent flow and calculate its slow entropy.

Let $G = \operatorname{SL}(d,\mathbb{R})$, $\Gamma \subset G$ be any lattice, and
	
	\[ U = \begin{pmatrix}
	0 & 1 \\
	& 0 &  1 \\
	& & \ddots & \ddots \\
	& & & 0 & 1 \\
	& & & & 0
	\end{pmatrix}, \qquad \exp_{alg}(tU) = \begin{pmatrix}
	1 & t & t^2/2 & \cdots & t^{d-1}/(d-1)! \\
	& 1 & t & \cdots & t^{d-2}/(d-2)! \\
	& & \ddots & \ddots & \vdots \\
	& & & 1 & t \\
	& & & & 1
	\end{pmatrix}.\]
	
	We call $U$ the {\it principal} nilpotent element associated to the algebra $\mathfrak{sl}(d,\mathbb{R})$. In the special case of $G=\operatorname{SL}(d,\mathbb{R})$, any nilpotent algebra element is conjugated to a block form element
	\[ U = \begin{pmatrix}
	U_1 \\
	& U_2 \\
	& & \ddots \\
	& & & U_n
	\end{pmatrix},
	\]
where each $U_i \in \mathfrak{sl}(d_i,\mathbb{R})$ is a principal element. Note that this is exactly the Jordan normal form of the matrix $U$. Call each $U_i$ a {\it principal block} and the sequence $(\dim(U_1),\dots,\dim(U_n))$ the {\it block sequence of $U$}.

We get the following Corollary of Theorem \ref{thm:main}:

\begin{corollary}(Kanigowski-Vinhage-Wei \cite{KanigowskiVinhageWei})
\label{cor:specific}
	The topological and measure-theoretic slow entropy of a principal unipotent flow equals $GR(U)$, where $$GR(U)=\sum_{i=1}^{m}\frac{1}{6}k_i(4k_i+1)(k_i-1)+\sum_{i=1}^{m-1}\sum_{j=i+1}^m\frac{1}{3} k_i \left(k_i^2+3k_j^2-3k_j-1\right), $$ where $(k_1,\dots,k_m)$, $k_i \le k_{i+1}$ is the block sequence of $U$.\label{sld-formula}
\end{corollary}

\subsection{Slow entropy of rank one systems}
Rank one systems form a very important class of systems in ergodic theory; they have a quite simple constructive definition and they are a rich source of examples and counterexamples in ergodic theory. We define a rank one system as follows:

\begin{definition}[$\beta$-Rank one]
Let $\beta$ belong to $(0,1]$. A system $(X,T,\mu)$ is of $\beta$-rank one if for every partition $\alpha$ of $X$, for every positive $\epsilon$, there exist two subsets $F$ and $A$ of $X$, a positive integer $h$ and a partition $\alpha'$ of $A$ such that,
\begin{enumerate}[(1)]
\item $F,TF,\ldots,T^{h-1}F$ are disjoint;
\item $A=\cup_{j=0}^{h-1}T^jF$;
\item $\mu(A)>\beta$;
\item $|\alpha|_A-\alpha'|<\epsilon$,
\item $\alpha'$ is refined by the sets $\{F,TF,\ldots,T^{h-1}F\}$.
\end{enumerate}
Moreover, if $\beta=1$, we will say $\beta-$rank one system is a rank one system.
\end{definition}

Recall that rank one systems have zero entropy, thus it is natural to use slow entropy to measure the complexity of rank one systems. By Ferenczi's results, we have the following general theorem about the upper bound of slow entropy for rank one systems:

\bt[Ferenczi \cite{Ferenczi}]\label{thm:slowEntRankone1}
For any rank one system $(X,T,\mu)$ and $a_n(t)=n^t$, we have $$\underline{\operatorname{ent}}_a^{\mu}(T)\leq2.$$
\et

Kanigowski \cite{Kanigowski} got more precise upper bounds on slow entropy of local rank one flows\footnote{Slow entropy can be defined analogously for $\mathbb{R}$-actions with the use of continuous Hamming metric.}:
\bt[Kanigowski \cite{Kanigowski}]\label{thm:rankoneSlowEntropy}
Let $\beta\in(0,1]$ and $g(n)$ be any sequence of positive numbers such that $\lim_{n\to+\infty}g(n)=+\infty$. Then for any measure-preserving $\beta-$rank one flow $T_t$ which acts on $(X,\mu)$ and $a_n(t)=n(g(n))^t$, we have $\underline{ent}_{a}^{\mu}(T_t)=0$.
\et

It is worth to point out that although the lower \textit{symbolic complexity} can be controlled for a rank one system with the same bound as above, the upper symbolic complexity with respect to the scaling function $a_n(t)=n^t$ can be $+\infty$. For more details, interested readers may go to \cite{Ferenczi96}.

For a special rank one system, namely the Chacon system, the precise value of slow entropy is obtained:

\begin{proposition}[Ferenczi \cite{Ferenczi}]
For the Chacon system $(X,T,\mu)$ and $a_n(t)=n^t$,
\begin{equation}
\underline{ent}_a^{\mu}(T)=\overline{ent}_a^{\mu}(T)=1,
\end{equation}
\end{proposition}

\subsection{Slow entropy of surface flows}
Smooth flows on surfaces play a really important role in dynamical systems. One of the reasons is that dimension $2$ is the lowest dimension in which we can have some non-trivial ergodic behavior. We restrict our attention to smooth surface flow with fixed points \footnote{If a smooth surface flow has no fixed points, then by Lefschetz formula, the surface is a two-dimensional torus and the flow is a smooth time-change of a linear flow.}. Recall that smooth surface flows have zero entropy due to Pesin formula, thus it is natural to try to compute slow entropy in this class. It turns out that ergodic properties of smooth surface flows are successfully  studied via their \textit{special representation}. We will restrict our attention to so called  Arnol'd and Kochergin type special flows (defined below) over irrational rotations. Kanigowski \cite{Kanigowski}, computed the precise value of slow entropy of these flows in the scales $a_n(t)=n(\log n)^t$ (for Arnol'd flows) and $a_n(t)=n^t$ (for Kochergin flows). Flows that we consider have the following special representation (with the base automorphism $T$ and the roof function $f$):

\bd[Arnol'd flows and Kochergin flows]
Let
\begin{itemize}
    \item $T=R_\alpha:\mathbb{T} \to\mathbb{T}$,  $R_\alpha x=x+\alpha \mod 1$;
	 \item $f$ is a $C^2(\mathbb{T}\setminus\{0\})$ function which satisfies:
\begin{equation}\label{Kor1}
\lim_{x\to0^+}\frac{f(x)}{h(x)}=A_1\text{ and }\lim_{x\to0^-}\frac{f(x)}{h(1-x)}=B_1, \text{ where }A_1,B_1>0;
\end{equation}
\begin{equation}\label{Kor2}
\lim_{x\to0^+}\frac{f'(x)}{h(x)}=-A_2\text{ and }\lim_{x\to0^-}\frac{f'(x)}{h(1-x)}=B_2, \text{ where }A_2,B_2>0;
\end{equation}
\begin{equation}\label{Kor3}
\lim_{x\to0^+}\frac{f''(x)}{h(x)}=A_3\text{ and }\lim_{x\to0^-}\frac{f''(x)}{h(1-x)}=B_3, \text{ where }A_3,B_3>0.
\end{equation}
\end{itemize}
\begin{itemize}
\item We call such flows {\em Kochergin flows} if $f(x)=x^{-\gamma}$, $0<\gamma<1$ and denote them by $\mathscr{T}^{f,\gamma}_t$.
\item We call such flows {\em Arnol'd flow} if $f(x)=-\log x$ and we assume additionally that $A_i+B_i\neq0, i=1,2,3$ and denote them as $\mathscr{T}_t^f$.
\end{itemize}
\ed

The following two full measure sets are needed for the main results:

\bd
Let $\{q_n\}$ be the denominators of $\alpha$ and let
$$\mathscr{D}:=\{\alpha\in\mathbb{R}\backslash\mathbb{Q}:q_{n+1}\leq C(\alpha)q_n\log q_n(\log n)^2\text{ for every } n\in \mathbb{N}\},$$
$$K_{\alpha}:=\{n:q_{n+1}\leq q_n\log^{\frac{7}{8}}q_n \text{ for every } n\in \mathbb{N}\},$$
and
$$\mathscr{E}:=\{\alpha\in\mathbb{R}\backslash\mathbb{Q}:\sum_{i\notin K_{\alpha}}\log^{-\frac{7}{8}}q_i<+\infty\}.$$
\ed
\brm
By Khinchin theorem \cite{Khinchin} we have $\lambda(\mathscr{D})=1$. By Fayad and Kanigowski \cite{KanigowshiFayad}, we have $\lambda(\mathscr{E})=1$.
\erm
Then we have following two theorems:
\bt[Kanigowski \cite{Kanigowski}]
Let $a_n(t)=n(\log n)^t$. Then for every $\alpha\in\mathscr{D}\cap\mathscr{E}$ and the corresponding Arnol'd flow $\mathscr{T}_t^f$, we have $\underline{ent}_a^{\mu}(\mathscr{T}_t^f)=1$.
\et

\bt[Kanigowski \cite{Kanigowski}]
Let $a_n(t)=n^t$. Then for every $\alpha\in\mathscr{D}$ and every $\gamma\in(0,1)$ for the corresponding Kochergin flow $\mathscr{T}^{f,\gamma}_t$, we have $\underline{ent}_a^{\mu}(\mathscr{T}^{f,\gamma}_t)=1+|\gamma|$.
\et

\subsection{Slow entropy of AbC constructions}\label{sec:slowEntropyABC}
Recall in Section \ref{sec:failureVariational} that  the AbC construction helps us to show the failure of variational principle for slow entropy.
\subsubsection{Brief description of the AbC construction}
We give a brief description of the Anosov-Katok construction on a two dimensional torus $\mathbb{T}^2=\mathbb{R}^2/\mathbb{Z}^2$. Let $S_t(x,y)=(x+t,y)$.

We use notation
\begin{equation}
    \Delta_q=\{(x,y)\in\mathbb{T}^2: 0\leq x\leq\frac{1}{q}\},\quad \eta_q=\{S_\frac{i}{q}\Delta_q\}_{i=0}^{q-1}.
\end{equation}

Given any summable sequence of positive real numbers $\{\epsilon_n\}$, the construction proceeds inductively. Assume that we have chosen sequences of integers $\{k_m\}$ and $\{l_m\}$ for $m=1,2,3,\ldots, n-1$, sequences of integers $\{p_m\}, \{q_m\}$, a sequence of rationals $\{\alpha_m\}$, sequences of diffeomorphisms $\{h_m\}, \{H_m\},\{T_m\}$, sequences of partitions $\{\eta_{m}\},\{\xi_{q_{m}}\}$, and a sequence of sets $\{E_m\}$ for $m=1,2,3,\ldots, n$, such that the following properties are satisfied for $m<n$:
\begin{equation}
\begin{aligned}
    & \alpha_m=\frac{p_m}{q_m},\quad \alpha_{m+1}=\alpha_m+\beta_m,\quad \beta_m=\frac{1}{k_ml_mq_m^2},\\
    & p_{m+1}=k_ml_mp_m+1,\quad q_{m+1}=k_ml_mq_m^2,\\
    & H_{m+1}=h_{m+1}\circ h_m\circ \ldots \circ h_1,\\
    & T_{m+1}=H_m^{-1}\circ S_{\alpha_m}\circ H_m,\\
    & \xi_{m+1}=H_{m+1}^{-1}\eta_{q_{m+1}},\\
    & \mu(E_m)<\epsilon_m,\\
    & \|T_{m+1}-T_m\|_{C^m}<\epsilon_m.
\end{aligned}
\end{equation}
At the $n+1$-th step we choose $k_n, l_n$,  $h_{n+1}$ and $E_{n+1}$. Note that this automatically gives us $q_{n+1}$, $p_{n+1}$ $H_{n+1}$, $\eta_{q_{n+1}}$, $\xi_{q_{n+1}}$, $\beta_{n+1}$, $\alpha_{n+1}$ and $T_{n+1}$. We finally define $T:=\lim_{n\to\infty}T_n$.

We have the following result regarding the slow entropy of AbC constructions:

\begin{theorem}[Kanigowski-Vinhage-Wei, \cite{KanigowskiVinhageWei3}]\label{thm:ABCentropy}
Assume that $\epsilon_{n_0}> \sum_{n=n_0+1}^\infty\epsilon_n$ for all $n_0$ sufficiently large. Then for any scale $a_n(t)$, there exists a diffeomorphism obtained by the AbC construction whose slow entropy is zero.
\end{theorem}

\brm
Here the topological slow entropy cannot vanish at all scales as (non-trivial) AbC constructions are not topologically conjugate to a translation.
\erm

Combining slow entropy Goodwyn's theorem (Theorem \ref{thm:Goodwyn}) with the above theorem, we have the following corollary.

\bc
Measure-theoretic slow entropy of an AbC construction can be made zero at arbitrary small scales by adjusting the coefficients of AbC construction.
\ec




\subsection{Slow entropy of higher rank actions}\label{sec:slowEntropyHighRank}

In order to measure the complexity of smooth measure-preserving (higher rank) actions, we consider the slow entropy invariant.  One needs to modify slow entropy original definition to a specific scale which will work better for higher rank abelian action:
\bd[Slow entropy of a $\mathbb{Z}^k$ action \cite{KatokKatokRH}]
Let $\alpha$ be an action of $\mathbb{Z}^k$. Given a norm $d$ on $\mathbb{R}^k$ let $F_s^d$ be the intersection of the lattice $\mathbb{Z}^k$ with the ball centered at $0$ and radius $s$. We define the slow entropy of $\alpha$ w.r.t. to the norm $d$ and the partition $\xi$ as

$$sh(d,\alpha,\xi)=\lim_{\epsilon,\delta\to0}\limsup_{s\to\infty}\frac{1}{s}\log S_{\xi}^H(\alpha,F_s^d,\epsilon,\delta).$$
By taking the supremum  over all finite measurable partitions, we define
\beq
sh(d,\alpha)=\sup_{\xi}sh(d,\alpha,\xi).
\ee
Finally, we define the slow entropy of $\alpha$ as
$$sh(\alpha)=\inf_{d:vol(d)=1}sh(d,\alpha),$$
where $vol(d)$ is the volume of the unit ball in the norm $d$ under the assumption that $\mathbb{Z}^k\subset\mathbb{R}^k$ has co-volume $1$.
\ed

In order to state results related to slow entropy of Cartan actions, we need to introduce some basic notation related to algebraic number theory.
Suppose $A\in\operatorname{SL}(n,\mathbb{Z})$ is a matrix and has an irreducible characteristic polynomial $f$. Denote the centralizer of $A$ with rational coefficients as $Z_Q(A)$, then $Z_Q(A)$ can be identified with the field $K=\mathbb{Q}(\lambda)$, where $\lambda$ is an eigenvalue of $A$. Denote this identification by $\gamma:p(A)\to p(\lambda)$ with $p\in\mathbb{Q}(x)$. Then define $O_K$ as the ring of integers in $K$ and $U_K$ as the group of units in $O_K$. We state the result on slow entropy of Cartan actions:
\bt[A.Katok-S.Katok-Rodriguez Hertz  \cite{KatokKatokRH}]
Let $\alpha$ be a Cartan action on the torus $\mathbb{T}^n$. Then
\beq
sh(\alpha)=C(n)R^{\frac{1}{n-1}},
\ee
with $R=kR_K$, where $R_K$ is the regulator, $k=[U_K:\gamma(\alpha)]$, and
\beq
\frac{n-1}{2}\leq C(n)\leq n-1.
\ee
\et

\section{Other entropy type invariant}\label{sec:Otherinvariants}
In this section, we will introduce other invariants related to slow entropy: Fried entropy, entropy convergence rate, Kakutani invariant and entropy dimension.
\subsection{Fried average entropy}
In this section we will define Fried average entropy for higher rank actions and then introduce some relations between slow entropy and Fried average entropy. We will follow \cite{Fried} and \cite{KatokKatokRH} in this section. We define the entropy function of a one parameter subgroup of $\mathbb{R}^k$ action which  preserves a probability ergodic measure $\mu$ at first:
\bd
The entropy function of an $\mathbb{R}^k$- action $\alpha$, denoted by $h_{\mu}^{\alpha}$, associates to $t^v=(t_1,\ldots,t_k)\in\mathbb{R}^k$ the value of the measure-theoretic entropy of $\alpha(t^v)$: $h_{\mu}^{\alpha}(t^v)=h_{\mu}(\alpha(t^v))$.
\ed

We now define the Fried average entropy of a $\mathbb{R}^k$ action
\bd[The Fried average entropy, Fried \cite{Fried}]
Let $\alpha$ be an $\mathbb{R}^k$ action with a fixed volume element. The {\em Fried average entropy} of $\alpha$, denoted as $h^F(\alpha)$, is defined as
$$h^F(\alpha)=\frac{2^k}{k!vol(B(h_{\mu}^{\alpha}))},$$
where $B(h_{\mu}^{\alpha})$ is the unit ball with respect to the entropy function.

If $\alpha$ is a $\mathbb{Z}^k$ action, then the Fried average entropy of $\alpha$ is defined as the Fried average entropy of its suspension.
\ed

Here are some basic properties of Fried average entropy:
\bp[\cite{Fried},\cite{KatokKatokRH}]
\begin{enumerate}[(a)]
\item Let $\alpha_1$ be an action of $\mathbb{R}^{k_1}$ on $M_1$, $\alpha_2$ be an action of $\mathbb{R}^{k_2}$ on $M_2$, preserving the measure $\mu_1$ and $\mu_2$ respectively, and $\alpha=\alpha_1\times\alpha_2$ be the product action of $\mathbb{R}^{k_1+k_2}$ with the measure $\mu=\mu_1\times\mu_2$. Then $$h^F(\alpha)=h^F(\alpha_1)\cdot h^F(\alpha_2).$$

\item For $k=1$, $h^F(\alpha)$ is equal to the usual entropy of the flow, i.e. the entropy of its time-one map.

\item The Fried average entropy is equal to zero if the entropy function is a semi-norm.

\item The Fried average entropy of a $\mathbb{Z}^k$ action is independent of the choice of generators.

\item The Fried average entropy of the restriction of a $\mathbb{Z}^k$ action $\alpha$ to a finite index subgroup $A\subset\mathbb{Z}^k$ is equal to $h_{F}^{\alpha}$ multiplied by the index of the subgroup $A$.
\end{enumerate}
\epo

Recall that $\alpha$ is a maximal rank action, if it is an action of  $\mathbb{Z}^a\times\mathbb{R}^b$ with $a+b=n-1$ on a $n+b-$dimensional manifold. We  have the following dichotomy for $k\geq 2$:

\bt[A.Katok-S.Katok-Rodriguez Hertz  \cite{KatokKatokRH}]
Suppose $\alpha$ is a maximal rank $\mathbb{Z}^k$ action preserving $\mu$. Then either $\mu$ is absolutely continuous and the Fried average entropy is positive or Fried average entropy of the action is equal to zero.
\et

For Cartan actions, we have the following quantitative description of Fried average entropy, we define Cartan action at first.
\bd[Cartan action, Hurder \cite{Hurder}]
Let $\mathscr{A}$ be a free abelian group with a given set of generators $\Delta=\{\gamma_1,\ldots,\gamma_n\}$. Then $(\phi,\Delta)$ is a Cartan $C^r$-action on the $n$-manifold $X$ if
\begin{enumerate}[(a)]
  \item $\phi:\mathscr{A}\times X\to X$ is a $C^r$-action on $X$;
  \item each $\gamma_i\in\Delta$ is $\phi$-hyperbolic, i.e. $\phi(\gamma_i)$ is an Anosov diffeomorphism of $X$, and $\phi(\gamma_i)$ has one-dimensional, strong stable foliation $\mathscr{F}_i^{ss}$;
  \item the tangential distributions $E_i^{ss}=T\mathscr{F}_i^{ss}$ are pairwise transverse with their internal direct sum $E_1^{ss}\oplus\ldots\oplus E_N^{ss}\cong TX$
\end{enumerate}
Moreover, if $\phi:\Gamma\times X\to X$ is an Anosov $C^r$-action on a manifold $X$, we say that $\phi$ is a Cartan lattice action if there is a subset of commuting hyperbolic elements $\Delta=\{\gamma_1,\ldots,\gamma_n\}\subset \Gamma$, which generate an abelian subgroup $\mathscr{A}$, such that the restriction of $\phi|\mathscr{A}$ is an abelian Cartan $C^r$-action on $X$.
\ed
Then we have,
\bt[Lower bound of Fried average entropy, A.Katok-S.Katok-Rodriguez Hertz \cite{KatokKatokRH}]
The Fried average entropy $h^F(\alpha)$ of a Cartan action $\alpha$ of a given rank $n-1\geq 2$ is bounded away from zero by a positive function that grows exponentially with $n$,
\beq\label{eq:exponentialGrowth}
h^F(\alpha)>0.000752\exp(0.244n),
\ee
moreover,
$$h^F(\alpha)\geq0.089.$$
\et

Combining the main result in \cite{KatRH} with the above theorem, we have the following characterization of weakly mixing maximal rank $\mathbb{Z}^k$ actions:
\bc[A.Katok-S.Katok-Rodriguez Hertz  \cite{KatokKatokRH}]
The Fried average entropy of a weakly mixing maximal rank action of $\mathbb{Z}^k$, $k\geq2$, is either equal to zero or is greater than $0.089$; furthermore the lower bound grows exponentially with $k$ as \eqref{eq:exponentialGrowth}.
\ec

At the end of this section, we will establish the connection between Fried average entropy and slow entropy of higher rank actions. The first relation is for affine actions on the torus:
\bt[A.Katok-S.Katok-Rodriguez Hertz  \cite{KatokKatokRH}]
For an affine action $\alpha$ on a torus, $sh(\alpha)=0$ if and only if $h^F(\alpha)=0$.
\et

If we assume $\alpha$ is a weakly mixing maximal rank action, we  get the following formulas explicitly established the relation between slow entropy and Fried average entropy:
\bt[A.Katok-S.Katok-Rodriguez Hertz  \cite{KatokKatokRH}]\label{thm:slowFriedConnection}
Suppose $\alpha$ is a weakly mixing maximal rank $\mathbb{Z}^k$ action for $k\geq2$. Then we have
\beq
sh(\alpha)=c(k)\binom{2k}{k}^{\frac{1}{k}}(h^F(\alpha))^{\frac{1}{k}},
\ee
where $c(k)\in[\frac{k}{4},\frac{k}{2}]$.
\et
\brm
Theorem \ref{thm:slowFriedConnection} implies that the slow entropy is either equal to zero or uniformly bounded away from zero. Moreover, the lower bounds grows linearly with respect to the rank.
\erm

\subsection{Kakutani equivalence invariant}\label{KakutaniEquivalenceInvariant}
Kakutani equivalence relation was first introduced by Kakutani \cite{Kakutani} in $1943$. Later this equivalence relation was developed by Katok \cite{KatokTimeChange}, \cite{KatokMonotone} and Feldman \cite{Feldman} by introducing a property named standardness (zero entropy loosely Bernoulli or loosely Kronecker) which is invariant under Kakutani equivalence.

Later, Ratner introduced a Kakutani invariant in \cite{RatnerKakInv} and \cite{RatnerKakInv2} to study Kakutani equivalence quantatively. In this section we will concentrate on this invariant and related results. The general idea of this invariant is to use $\bar{f}-$metric instead of Hamming metric to measure the complexity of a given system.


\subsubsection{Definition and basic theorems}
Let $T$ and $S$ be two ergodic transformations on $(X,\mathscr{B},\mu)$ and $(Y,\mathscr{C},\nu)$, respectively. We call $T$ and $S$ \emph{Kakutani equivalent}, denoted as $T\overset{K}{\sim}S$, if there exist $A\subset X$ and $B\subset$ with $\mu(A)\nu(B)>0$ such that $(T|_A,A,\mu_A)$ and $(S|_B,B,\nu_B)$ are measurably isomorphic.

Kakutani equivalence can also be defined through \emph{time changes}. For a flow $T_t$ acting on $(X,\mathscr{B},\mu)$ and a function $\alpha\in L^1_+(X,\mathscr{B},\mu)$, the flow $T_t^{\alpha}$ is called a \emph{time change} of $T_t$ if $T_t^{\alpha}(x)=T_{u(t,x)}(x)$, where $u(t,x)$ satisfies $\int_0^{u(t,x)}\alpha(T_sx)ds=t$. Then two ergodic measure preserving flows $T_t$, $S_t$ acting on $(X,\mathscr{B},\mu)$ and $(Y,\mathscr{C},\nu)$ are \emph{Kakutani equivalent} if $S_t$ is isomorphic to $T_t^{\alpha}$ for some $\alpha\in L_+^1(X,\mathscr{B},\mu)$.

In $1977$, Katok introduced a hierarchy structure among Kakutani equivalence classes by a binary relation called \emph{majorization}. An ergodic flow $S_t$ is majorized by an ergodic flow $T_t$, which is denoted as $S_t\overset{M}{\prec}T_t$, if there exists an ergodic flow $\tilde{T}_t$ which is Kakutani equivalent to $T_t$ such that $S_t$ is a factor of $\tilde{T}_t$. It is clear that if $T_t\overset{K}{\sim}S_t$, then we have $S_t\overset{M}{\prec}T_t$ and $T_t\overset{M}{\prec}S_t$, which guarantees that majorization defines a quasi-partial order among all Kakutani equivalence classes\footnote{We do not know whether $S_t\overset{M}{\prec}T_t$ and $T_t\overset{M}{\prec}S_t$ imply $S_t\overset{K}{\sim}T_t$.}.

Among all Kakutani equivalence classes, there exists a special equivalence class majorized by all other classes which is called the \emph{standard class} and can be characterized geometrically as follows: a flow $T_t$ is {\em standard} or {\it loosely Bernoulli of zero entropy} if it is Kakutani equivalent to an irrational linear flow \footnote{Notice that standardness implies ergodicity. By Katok's results in \cite{KatokMonotone}, all irrational linear flows on $\mathbb{T}^2$ are Kakutani equivalent.} on $\mathbb{T}^2$.

From now on, we will assume that $T_t$ is an ergodic flow acting on a Lebesgue space $(X,\mathscr{B},\mu)$, $\mathcal{P}$ is a finite measurable partition of $X$ and $I_T(x)=\{T_sx:s\in[0,T]\}$. Suppose that $x\in X$ and we denote $\mathcal{P}(x)$ as the atom of $\mathcal{P}$ containing $x$.

\begin{definition}[$(\epsilon,\mathcal{P})$-matchable, Ratner \cite{RatnerKakInv2}]
For $x,y\in X$, $\epsilon>0$ and $T>1$, $I_T(x)$ and $I_T(y)$ are called $(\epsilon,\mathcal{P})$-matchable if there exist subsets $A=A(x,y),A'=A'(x,y)\subset[0,T]$ with $l(A)>(1-\epsilon)T$, $l(A')>(1-\epsilon)T$ and an increasing absolutely continuous map $h=h(x,y)$ from $A$ onto $A'$ such that $\mathcal{P}(T_tx)=\mathcal{P}(T_{h(t)}y)$ for all $t\in A$ and the derivative $h'$ satisfies
\beq
|h'(t)-1|<\epsilon,\ \ \forall t\in A.
\ee
We call $h$ an $(\epsilon,\mathcal{P})-$matching from $I_T(x)$ onto $I_T(y)$.
\end{definition}

With the above definition, Kakutani invariants can be defined in the following process (following Ratner):

\begin{definition}[Kakutani invariant, Ratner \cite{RatnerKakInv2}]
We have:
\begin{enumerate}[(1)]
\item Define \begin{equation}\label{eq:fbarContinuous}
    \bar{f}_T(x,y,\mathcal{P})=\inf\{\epsilon:\text{$I_T(x)$ and $I_T(y)$ are $(\epsilon, \mathcal{P})$-matchable}\};
    \end{equation}
\item Denote $B_T(x,\epsilon,\mathcal{P})=\{y\in X:\bar{f}_T(x,y,\mathcal{P})<\epsilon\}$ the $(T,\mathcal{P})$-ball of radius $\epsilon>0$ centered at $x\in X$, $T>1$;
\item A family $\alpha_T(\epsilon,\mathcal{P})$ of $(T,\mathcal{P})$-balls of radius $\epsilon>0$ is called $(\epsilon,T,\mathcal{P})$-cover of $X$ if we have $\mu(\cup\alpha_T(\epsilon,\mathcal{P}))>1-\epsilon$;
\item Denote $K_T(\epsilon,\mathcal{P})=\inf|\alpha_T(\epsilon,\mathcal{P})|$ where $|A|$ denote the number of elements in $A$ and $\inf$ is taken over all $(\epsilon, T, \mathcal{P})$-covers of $X$.
\item Let $\mathcal{F}$ denote the family of all nondecreasing functions from $\mathbb{R}^+$ onto itself, converging to $+\infty$.
\item For $u\in \mathcal{F}$, we denote,
\begin{equation}
\begin{aligned}
\beta(u,\epsilon,\mathcal{P})&=\liminf_{T\to\infty}\frac{\log K_T(\epsilon, \mathcal{P})}{u(T)};\\
e(u,\mathcal{P})&=\limsup_{\epsilon\to0}\beta(u,\epsilon,\mathcal{P});\\
e(T_t,u)&=\sup_{\mathcal{P}}e(u,\mathcal{P}).
\end{aligned}
\end{equation}
\end{enumerate}
\end{definition}

We shortly document some basic properties of $e(T_t,u)$ in the following:
\bt[Generator theorem, Ratner \cite{RatnerKakInv}, \cite{RatnerKakInv2}]
Let $T_t$ be an ergodic measure-preserving flow on $(X,\mathscr{B},\mu)$ and let $\mathcal{P}_1\leq\mathcal{P}_2\leq\ldots$ be an increasing sequence of finite measurable partitions of $X$ such that $\vee_{n=1}^{\infty}\mathcal{P}_n$ generates the $\sigma-$algebra $\mathscr{B}$. Then $e(T_t,u)=\sup_m e(u, \mathcal{P}_m)$ for all $u\in\mathcal{F}$.
\et


\bt[Ratner \cite{RatnerKakInv}, \cite{RatnerKakInv2}]\label{thm:KakutaniInvariantScaleFunction}
Let $T_t$ and $\tilde{T}_t$ be two ergodic Kakutani equivalent measure preserving flows on $(X,\mathscr{B},\mu)$ and $(\tilde{X},\tilde{\mathscr{B}},\tilde{\mu})$ respectively. Then
$$e(T_t,u)=e(\tilde{T}_t,u)$$
for all $u\in U$ with
$$\lim_{t\to\infty}\frac{u(at)}{u(t)}=1\text{ for all }a>0.$$
\et
\brm
Theorem \ref{thm:KakutaniInvariantScaleFunction} guarantees $e(T_t,u)$ is a Kakutani invariant for some certain functions $u$. A simple example which satisfies the assumptions of Theorem \ref{thm:KakutaniInvariantScaleFunction} condition is $u(t)=\log t$.
\erm

\subsubsection{Application of Kakutani invariants}

With the help of Kakutani invariant, we can give another characterization of standard flows, which should be understood as a quantitative version of majorization characterization:
\bt[Ratner \cite{RatnerKakInv}, \cite{RatnerKakInv2}]\label{th:LB}
A zero-entropy ergodic measure preserving flow $T_t$ is standard if and only if $e(T_t,u)=0$ for all $u\in \mathcal{F}$.
\et

Recall that Ratner proved that the horocycle flow is standard and its Cartesian product is not standard in \cite{RatnerHoroLB} and \cite{RatnerHoroPrNLB}. Moreover, Ratner gave the following estimates of the Kakutani invariant of $n-$fold Cartesian product of the horocycle flow.

\bt[Ratner \cite{RatnerKakInv}, \cite{RatnerKakInv2}]\label{thm:nhorocycleflow}
Let $h_t^{(n)}=h_t\times\cdots\times h_t$ $n=1,2,\ldots,$ be the $n-$fold Cartesian product of the horocycle flow on $\left(\operatorname{SL}(2,\mathbb{R})/\Gamma\right)^n$ with cocompact lattice $\Gamma$ and $u(t)=\log t$, $t>1$. Then
$$3n-3\leq e(h^{(n)}_t,u)\leq 3n-2,$$
for all $n=1,2,\ldots$.
\et

Combining Theorem \ref{th:LB} and Theorem \ref{thm:nhorocycleflow}, it is clear that among the $n-$fold Cartesian product of the horocycle flow, only the horocycle flow itself is standard.

Inspired by M.\ Ratner's results on the horocycle flow, it is natural to consider the Kakutani equivalence for general unipotent flows, this has even been considered by M. Ratner  in  \cite{RatICM}. Based on multiscale analysis and structure of semisimple Lie groups, Kanigowski, Vinhage and Wei obtained the following result about the Kakutani invariant for unipotent flows. Recall that $GR(U)$ is defined in Theorem \ref{thm:slowentropyinvariants}.

\begin{theorem}[Kanigowski-Vinhage-Wei \cite{KanigowskiVinhageWei2}]\label{thm:KakutaniInvariantsUnipotents}
Let $G$ be a semisimple linear Lie group and
$\phi_t= L_{\exp(tU)}$ a unipotent flow on $G\slash \Gamma$. If $\Gamma$ is cocompact, we have
$$
e(\phi_t,\log)=GR(U)-3.
$$
For finite volume $\Gamma$, we have
$$GR(U)-4\leq e(\phi_t,\log)\leq GR(U)-3.$$
Moreover, if $GR(U)=3$, then $(\phi_t)$ is standard.
\end{theorem}

By a direct computation, we get $3n-4\leq e(h_t^{(n)},\log)\leq 3n-3$. This generalises M. Ratner's results \cite{RatnerKakInv} and \cite{RatnerKakInv2} to any lattice in $\left(\operatorname{SL}(2,\mathbb{R})\right)^k$. If we assume additionally that the lattice is cocompact, then $e(h_t^{(n)},\log)=3n-3$. Moreover, combining Theorem \ref{thm:KakutaniInvariantsUnipotents} with semisimple Lie group structure, we can use the following corollaries to provide a solution to M. Ratner's ICM Problem 1 in \cite{RatICM} (see also \cite{KleinShahStarkov}):

\begin{corollary}[Kanigowski-Vinhage-Wei \cite{KanigowskiVinhageWei2}]\label{cor1:KakutaniInvariantsUnipotents} The only ergodic unipotent flows on compact quotients of linear semisimple Lie groups which are standard are of the form $\phi_t=\begin{pmatrix}1&t\\0&1\end{pmatrix}\times \operatorname{id}$ acting on $(\operatorname{SL}(2,\mathbb{R})\times G')\slash \Gamma$,
where $\Gamma$ is irreducible.
\end{corollary}

\begin{corollary}[Kanigowski-Vinhage-Wei \cite{KanigowskiVinhageWei2}]\label{cor2:KakutaniInvariantsUnipotents}
Let $G$ be a linear semisimple Lie group with $\dim G>3$, and $G / \Gamma$ be a compact homogeneous space of $G$.

\begin{enumerate}[(a)]
\item There are ergodic unipotent flows on $G/\Gamma$ which are not standard.
\item If $G$ is simple, no unipotent flow on $G/\Gamma$ is standard.
\item If $G$ has real rank at least two, there are two unipotent flows on $G/\Gamma$ (which are not identity, but not necessarily ergodic) which are not Kakutani equivalent.
\item If $G \cong \operatorname{SL}(d,\mathbb{R})$, then there are at least $d-1$ flows on $G/\Gamma$ which are pairwise non-Kakutani equivalent.
\end{enumerate}
\end{corollary}

\subsection{Scaled Entropy}
Scaled entropy, which was introduced by Pesin and Zhao (\cite{ZhaoPesin}, \cite{ZhaoPesin2}), is an invariant defined through the Carath\'{e}odory approach. This invariant also measures the complexity of a given system. Although scaled entropy may be quite different from slow entropy, there exist some relations which can be used to give an intuitive understanding of scaled entropy based on slow entropy. In the remaining part of this section, we give a short introduction of scaled entropy based on \cite{ZhaoPesin}, \cite{ZhaoPesin2}.

\subsubsection{Scaled topological entropy}
Let $(X,d)$ be a compact metric space and $T:X\to X$ be a continuous transformation.
\begin{enumerate}[(i)]
  \item Let $\mathbf{a}=\{a_n\}_{n\geq1}$ be a sequence of positive numbers which we call a \emph{scaled sequence} if it is monotonically increasing to infinity;
  \item Given  $\epsilon>0$, a scaled sequence $\mathbf{a}=\{a_n\}_{n\geq1}$, a subset $Z\subset X$ and $\alpha,N>0$, we define
\begin{equation}\label{eq:scaledEntropyM}
M(Z,\alpha,N,\epsilon,\mathbf{a})=\inf\left\{\sum_i\exp(-\alpha a_{n_i}):\cup_iB_{n_i}(x_i,\epsilon)\supset Z,\ \ x_i\in X\text{ and }n_i\geq N\ \ \forall i\right\},
\end{equation}
where $B_n(x_i,\epsilon)=\{y\in X:d_n^T(x,y)<\epsilon\}$ and $d_n^T$ is the Bowen metric defined in \eqref{eq:BowenMetric};
  \item As $M(Z,\alpha,N,\epsilon,\mathbf{a})$ is monotonically increasing with respect to $N$, its limit exists as $N\to\infty$ and we define $m(Z,\alpha,N,\epsilon,\mathbf{a})=\lim_{N\to\infty}M(Z,\alpha,N,\epsilon,\mathbf{a})$;
  \item Then define $E_Z(T,\mathbf{a},\epsilon)$ and $E_Z(T,\mathbf{a})$ as
\begin{equation}\label{eq:scaledTopoEntropy}
\begin{aligned}
E_Z(T,\mathbf{a},\epsilon)&=\inf\{\alpha:m(Z,\alpha,\epsilon,\mathbf{a})=0\}=\sup\{\alpha:m(Z,\alpha,\epsilon,\mathbf{a})=+\infty\};\\
E_Z(T,\mathbf{a})&=\lim_{\epsilon\to0}E_Z(T,\mathbf{a},\epsilon).
\end{aligned}
\end{equation}
We call $E_Z(T,\mathbf{a})$ the \emph{scaled topological entropy} of $T$ on the set $Z$.
  \item For $Z\subset X$, let $N(Z,n,\epsilon)$ be the smallest number of Bowen's ball $B_n(x,\epsilon)$ that needed to cover the set $Z$. Then for any scaled sequence $\mathbf{a}=\{a_n\}_{n\geq1}$, we define the following two quantities:
\begin{equation}\label{eq:UpLowScaledEntropyPre}
\begin{aligned}
\underline{E}_Z(T,\mathbf{a},\epsilon)&=\liminf_{n\to\infty}\frac{1}{a_n}\log N(Z,n,\epsilon),\\
\overline{E}_Z(T,\mathbf{a},\epsilon)&=\limsup_{n\to\infty}\frac{1}{a_n}\log N(Z,n,\epsilon).
\end{aligned}
\end{equation}
By taking $\epsilon\to0$, we define \emph{lower scaled topological entropy} $\underline{E}_Z(T,\mathbf{a})$ and \emph{upper scaled topological entropy} $\overline{E}_Z(T,\mathbf{a})$ as:
\begin{equation}\label{eq:UpLowScaledEntropy}
\begin{aligned}
\underline{E}_Z(T,\mathbf{a})&=\lim_{\epsilon\to0}\underline{E}_Z(T,\mathbf{a},\epsilon),\\
\overline{E}_Z(T,\mathbf{a})&=\lim_{\epsilon\to0}\overline{E}_Z(T,\mathbf{a},\epsilon).
\end{aligned}
\end{equation}
\end{enumerate}
\brm
Scaled entropy $E_Z(T,\mathbf{a})$, lower scaled topological entropy $\underline{E}_Z(T,\mathbf{a})$ and upper scaled topological entropy $\overline{E}_Z(T,\mathbf{a})$ are invariant under a topological conjugacy. Moreover, we have the following inequalities that describe the relation of these three invariants:
\beq
E_Z(T,\mathbf{a})\leq\underline{E}_Z(T,\mathbf{a})\leq\overline{E}_Z(T,\mathbf{a}).
\ee
\erm
\subsubsection{Scaled metric entropy}
If $\mu$ is a  $T-$invariant measure, we can define \emph{scaled metric entropy} for a given scaled sequence $\mathbf{a}$ as :
\beq\label{eq:scaledMetricEntropy}
E_{\mu}(T,\mathbf{a})=\lim_{\epsilon\to0}\lim_{\delta\to0}\inf\{E_Z(T,\mathbf{a},\epsilon):\mu(Z)>1-\delta\},
\ee
where $E_Z(T,\mathbf{a},\epsilon)$ is defined in
\eqref{eq:scaledTopoEntropy}.

Moreover, define \emph{lower scaled metric entropy} $\underline{E}_{\mu}(T,\mathbf{a})$ and \emph{upper scaled metric entropy} $\overline{E}_{\mu}(T,\mathbf{a})$ as
\beq\label{eq:UpLowScaledMetricEntropy}
\begin{aligned}
\underline{E}_{\mu}(T,\mathbf{a})&=\lim_{\epsilon\to0}\lim_{\delta\to0}\inf\{\underline{E}_Z(T,\mathbf{a},\epsilon):\mu(Z)>1-\delta\},\\
\overline{E}_{\mu}(T,\mathbf{a})&=\lim_{\epsilon\to0}\lim_{\delta\to0}\inf\{\overline{E}_Z(T,\mathbf{a},\epsilon):\mu(Z)>1-\delta\},
\end{aligned}
\ee
where $\underline{E}_Z(T,\mathbf{a},\epsilon)$ and $\overline{E}_Z(T,\mathbf{a},\epsilon)$ are defined in \eqref{eq:UpLowScaledEntropyPre}.

\brm
Same as in the topological case, scaled metric entropy, lower scaled metric entropy and upper scaled metric entropy are invariant under a topological  conjugacy. Moreover, we have
\beq
E_{\mu}(T,\mathbf{a})\leq\underline{E}_{\mu}(T,\mathbf{a})\leq\overline{E}_{\mu}(T,\mathbf{a}).
\ee
\erm
\subsubsection{Scaled entropy's properties}
By \cite{ZhaoPesin2}, we can describe the scaled metric entropy based on slow entropy language. Let $\mathbf{a}=\{a_n\}_{n\geq1}$ be a scaled sequence and $\mu$ is a $T-$invariant ergodic measure, we have
\beq
\begin{aligned}
\overline{S}(T,\mu)\leq\overline{E}_{\mu}(T,\mathbf{a}),\\
\underline{S}(T,\mu)\leq\underline{E}_{\mu}(T,\mathbf{a}),
\end{aligned}
\ee
where $\overline{S}(T,\mu)$ and $\underline{S}(T,\mu)$ defined as
\beq
\begin{aligned}
\overline{S}(T,\mu)=&\sup_{\xi}\lim_{\epsilon\to0}\lim_{\delta\to0}\limsup_{n\to\infty}\frac{\log S_{\xi}^H(T,n,\epsilon,\delta)}{a_n},\\
\underline{S}(T,\mu)=&\sup_{\xi}\lim_{\epsilon\to0}\lim_{\delta\to0}\liminf_{n\to\infty}\frac{\log S_{\xi}^H(T,n,\epsilon,\delta)}{a_n},
\end{aligned}
\ee
for $\xi$ from the collection of all finite measurable partitions and $S_{\xi}^H(T,n,\epsilon,\delta)$ from Section \ref{sec:SlowMeaureTHeoreticEntropy}.
\subsection{Entropy dimensions}\label{EntropyDimensions}
Similar as slow entropy, entropy dimension is also an invariant which is finer than the classical entrop. We give a basic introduction here and all material in this section is contained in  \cite{DeCarvalho}, \cite{FerencziPark}, \cite{ADP}, \cite{DouHuangPark} and \cite{DouHuangPark2}.

\bd[Topological entropy dimension]
Let $X$ be a compact metric space and $T$ be a continuous endormorphism of $X$. For an open cover $\mathscr{U}$ denote by $N(\mathscr{U})$ the minimal number of elements among all subcovers of $\mathscr{U}$ and $\mathscr{U}_0^{n-1}=\mathscr{U}\vee T^{-1}\mathscr{U}\ldots\vee T^{-n+1}\mathscr{U}$. Then the upper entropy dimension of $(X,T)$ w.r.t. an open cover is
\beq
\begin{aligned}
\overline{D}(\mathscr{U})&=\inf\{s\geq0:\limsup_{n\to\infty}\frac{1}{n^s}\log(N(\mathscr{U}_0^{n-1}))=0\}\\
&=\sup\{s\geq0:\limsup_{n\to\infty}\frac{1}{n^s}\log N(\mathscr{U}_0^{n-1})=\infty\}.
\end{aligned}
\ee
The upper entropy dimension of $(X,T)$ is defined as
$$\overline{D}(X,T)=\sup_{\mathscr{U}}\overline{D}(\mathscr{U}).$$

The lower entropy dimension of $(X,T)$ w.r.t. an open cover $\mathscr{U}$ is defined as
\beq
\begin{aligned}
\underline{D}(\mathscr{U})&=\inf\{s\geq0:\liminf_{n\to\infty}\frac{1}{n^s}\log(N(\mathscr{U}_0^{n-1}))=0\}\\
&=\sup\{s\geq0:\liminf_{n\to\infty}\frac{1}{n^s}\log N(\mathscr{U}_0^{n-1})=\infty\}.
\end{aligned}
\ee

The lower entropy dimension of $(X,T)$ is defined to be
$$\underline{D}(X,T)=\sup_{\mathscr{U}}\underline{D}(\mathscr{U}).$$

If $\overline{D}(X,T)=\underline{D}(X,T)$, we denote it as $D(X,T)$ and called the entropy dimension of $(X,T)$.
\ed

We also have measure-theoretic entropy dimension:
\bd[Measure-theoretic entropy dimension]
Let $(X,\mathscr{B},\mu,T)$ be a measure preserving system, $\alpha$ a finite measurable partition of $X$ and $\epsilon>0$. Denote $K(n,\epsilon)$ as the minimal cardinality of a set of Hamming balls\footnote{See \ref{def:HammingMetric} for more details.} with radius $\epsilon$ in $d_n^H$ distance that covers a set of measure at least $1-\epsilon$. Then define
$$\overline{D}(\alpha,\epsilon)=\sup\{s\in[0,1]:\limsup_{n\to\infty}\frac{\log K(n,\epsilon)}{n^s}>0\},$$
$$\overline{D}(\alpha)=\lim_{\epsilon\to0}\overline{D}(\alpha,\epsilon),$$
$$\overline{D}_{\mu}(X,T)=\sup_{\alpha}\overline{D}(\alpha).$$
and
$$\underline{D}(\alpha,\epsilon)=\sup\{s\in[0,1]:\liminf_{n\to\infty}\frac{\log K(n,\epsilon)}{n^s}>0\},$$
$$\underline{D}(\alpha)=\lim_{\epsilon\to0}\underline{D}(\alpha,\epsilon),$$
$$\underline{D}_{\mu}(X,T)=\sup_{\alpha}\underline{D}(\alpha).$$
$\overline{D}_{\mu}(X,T)$ and $\underline{D}_{\mu}(X,T)$ are called upper and lower entropy dimensions of $(X,\mathscr{B},\mu,T)$. We denote them as $D_{\mu}(X,T)$ if they are equal and calle as entropy dimension.
\ed

\brm
In \cite{DeCarvalho}, de Carvalho showed that if the topological entropy is finite then  both upper and lower entropy dimension is less than or equal to $1$ and if the topological entropy is finite and positive then both upper and lower entropy dimension is equal to $1$. From this result, we know that entropy dimension should be mainly used for zero entropy systems.
\erm

\brm
In \cite{DouHuangPark} and \cite{DouHuangPark2}, Dou, Huang and Park defined topological and measure-theoretic entropy dimension via the dimension of entropy generating sequences. Moreover, they also established the Goodwyn's inequality and many interesting properties of entropy dimension.
\erm


\brm
Ahn, Dou and Park \cite{ADP} proved that there exists a system such that topological entropy dimension of a topological dynamical system $(X,T)$ is strictly larger than the supremum of its measure-theoretic entropy dimensions, which demonstrates the failure of variational principle for entropy dimension.
\erm

Same as for slow entropy, entropy dimension also has product property in the non-classical version:
\bp[Ferenczi-Park \cite{FerencziPark}]
\beq
\overline{D}_{\mu}(X\times Y, T_1\times T_2)=\max\{\overline{D}_{\mu}(X, T_1),\overline{D}_{\mu}(Y,T_2)\}.
\ee
\epo

\brm
It is also worth to notice that measure-theoretic entropy dimension has a flexibility phenomenon. Namely for any $\alpha\in(0,1)$, there exists system $(X,T,\mathscr{B},\mu)$ such that $D_{\mu}(X,T)=\alpha$. For more details, we refer to \cite{FerencziPark} for reference.
\erm

\brm
It is also shown by Ferenczi and Park \cite{FerencziPark} that the upper and lower measure-theoretic entropy dimension may be different, i.e.\ there exists a system $(X,T,\mathscr{B},\mu)$ such that $\overline{D}_{\mu}(X,T)=1$ and $\underline{D}_{\mu}(X,T)=0$.
\erm

\subsection{Entropy convergence rate}\label{sec:EntropyConvergenceRate}
Instead of considering the Hamming  or Bowen metric, Blume (\cite{Blume1}, \cite{Blume2}, \cite{Blume3}) created a new isomorphism invariant, called \textit{entropy convergence rate}, for measure-preserving systems by using a direct approach to study  zero entropy systems, i.e.\ only changing the scaling function of original entropy function. The idea here is  straightforward: if there  is  no effective
upper sublinear bound among all partitions,  one looks at the lower bound instead. Unlike Kushnirenko's sequence entropy where vanishing of all invariants  characterizes systems with discrete spectrum,  Blume shows that there is a universal lower bound, somewhat slower that logarithmic for all  aperiodic transformations; for rank one systems it is at most  logarithmic and
for mixing rank one systems  the logarithmic lower bound is sharp. In order to avoid an obvious problem of missing more complicated
parts of the orbit structure, this approach needs to be modified; for example one may consider the slowest growth  rate of joint entropy only among generating partitions.

Let $(X,T,\mathscr{B},\mu)$ be a probability space with a measure-preserving endormorphism $T$ and $\alpha$ be a finite partition of $X$, then we define $\alpha^{n-1}=\vee_{i=1}^nT^{-i}\alpha$. By defining $f(x)=-x\log_2x$ for $x\in[0,1]$, the entropy of a partition $\alpha$ is defined as
$H(\alpha)=\sum_{A\in\alpha}f(\mu(A))$.
\bd[Entropy convergence rate, Blume \cite{Blume1}]
Let $\{a_n\}_{n\geq1}$ be a monotonically increasing sequence with $\lim_{n\to\infty}a_n=\infty$ and $c\in(0,\infty)$. If $\mathcal{P}$ is a class of partitions of $X$, then we say that $(X,T)$ is of type $(LS\geq c)$ for $(\{a_n\}_{n\geq1},\mathcal{P})$ if $$\limsup_{n\to\infty}\frac{H(\alpha^{n-1})}{a_n}\geq c,\text{ for all }\alpha\in \mathcal{P};$$
and $(X,T)$ is of type $(LI\geq c)$ for $(\{a_n\}_{n\geq1},\mathcal{P})$ if
$$\liminf_{n\to\infty}\frac{H(\alpha^{n-1})}{a_n}\geq c,\text{ for all }\alpha\in \mathcal{P}.$$
\ed
\brm
The types $(LS<\infty)$, $(LS=\infty)$, $(LI<\infty)$, $(LI=\infty)$, $(LI>0)$, $(LS\leq c)$, $(LI\leq c)$ are defined analogously.
\erm
\brm
By definition, we have following relations:
\beq
\begin{aligned}
&(LI\geq c)\Rightarrow(LS\geq c), \ \ (LS\leq c)\Rightarrow(LI\leq c),\\
&(LI>0)\Rightarrow(LS>0), \ \ (LS<\infty)\Rightarrow(LI<\infty).
\end{aligned}
\ee
\erm

In order to state results about the entropy convergence type of a measure preserving system, we need to specify the class of partitions under consideration. In the following text, we will mainly consider the following two cases:
\begin{enumerate}[(1)]
  \item $P(X)=\{\alpha(E)\;|\;\alpha(E)=\{E,X\setminus E\}\text{ for some }E\in\mathscr{B}\text{ with }0<\mu(E)<1\}$;
  \item $R(X)=\{\alpha(E)\;|\;E\in\mathscr{B}\text{ and }\lim_{n\to\infty}\max\{\mu(A)|A\in\alpha(E)^{n-1}_0\}=0\}$.
\end{enumerate}

Now we have enough preparation to introduce related results:
\bt[Blume  \cite{Blume1}]
If $(X,T)$ is an aperiodic measure-preserving system and $\{a_n\}_{n=1}^{\infty}$ is a positive monotone increasing sequence with $\lim_{n\to\infty}\frac{a_n}{n}=0$, then there exists an $E\subset\mathscr{B}$ such that
$$\lim_{n\to\infty}\frac{H(\alpha(E)^{n-1})}{a_n}=\infty.$$
\et

\brm
By definition, it is easy to see that this implies non-existence of aperiodic systems of types $(LS<\infty)$, $(LS\leq c)$ or $(LI\leq c)$.
\erm

However, if we assume additionally that $T$ is completely ergodic, then we can describe the entropy convergence rate for general monotonically increasing sequence $a_n$:
\bt[Blume \cite{Blume2}]
If $T$ is a completely ergodic transformation on $([0,1],\mathscr{B},\mu)$, then there exists a positive monotonically increasing sequence $\{a_k\}_{k\geq1}$ with $\lim_{k\to\infty}a_k=\infty$ such that $\forall E\in\mathscr{B}$ with $\mu(E)\in(0,1)$, we obtain
$$\liminf_{k\to\infty}\frac{H(\alpha(E)^{k-1})}{a_k}\geq1.$$
\et

On the other hand, if we keep aperiodicity condition but add some additional conditions on the sequence $a_n$ and the partition $\alpha(E)$, we will obtain almost same result:
\bt[Blume  \cite{Blume1}]\label{thm:entropyConvergenceRateAperiodic}
Let $(X,T)$ be an aperiodic measure-preserving system and let $g:[0,\infty)\to\mathbb{R}$ be a positive monotonely increasing function with $\int_1^{\infty}\frac{g(x)}{x^2}dx<\infty$. If $E\in\mathscr{B}$ is such that
$$\lim_{n\to\infty}\max\{\mu(A)|A\in\alpha(E)^{n-1}\}=0,$$
then
$$\limsup_{n\to\infty}\frac{H(\alpha(E)^{n-1})}{g(\log_2n)}=\infty.$$
\et

Moreover, if we replace aperiodicity by complete ergodicity, we will have the following characterization of completely ergodic systems based on the function $g$:
\bc[Blume \cite{Blume1}]\label{cor:entropyConvergenceRateComErgo}
Let $(X,T)$ be completely ergodic and let $g$ be given as in Theorem \ref{thm:entropyConvergenceRateAperiodic}. If $E\subset\mathscr{B}$ is such that $0<\mu(E)<1$, then
$$\limsup_{n\to\infty}\frac{H(\alpha(E)^{n-1})}{g(\log_2n)}=\infty.$$
\ec

Recall so far we always consider $g$ satisfying $\int_1^{\infty}\frac{g(x)}{x^2}dx<\infty$. If we consider some $g$ with some fast growth, i.e. $\int_1^{\infty}\frac{g(x)}{x^2}dx=\infty$, then we obtain:
\bt[Blume \cite{Blume2}]
If $T$ is a completely ergodic transformation on $([0,1],\mathscr{B},\mu)$, then there exists a positive concave function $g$ defined on $[1,\infty)$ with $\int_1^{\infty}\frac{g(x)}{x^2}dx=\infty$, such that for all $E\in\mathscr{B}$ with $\mu(E)\in(0,1)$, we obtain
$$\limsup_{k\to\infty}\frac{H(\alpha(E)^{k-1})}{g(\log_2k)}\geq1.$$
\et

Blume also used entropy convergence rate to characterize rank one systems:
\bt[Blume \cite{Blume1}]\label{thm:entropyConvergenceRateRankOne}
If $(X,T)$ is a rank-one system, then there is an $\alpha\in P(X)$ such that
$$\liminf_{n\to\infty}\frac{H(\alpha^{n-1})}{\log_2n}\leq2.$$

If we assume additionally that $(X,T)$ is mixing and $E\in\mathscr{B}$ with $0<\mu(E)<1$, then
$$\limsup_{n\to\infty}\frac{H(\alpha(E)^{n-1})}{\log_2n}>0.$$
\et

\brm
By Theorem \ref{thm:entropyConvergenceRateRankOne} it follows that no rank-one system can have a convergence rate of type $(LI\geq c)$ for any sequence $a_n$ that grows faster than $\log_2n$. Recall that Ferenczi (Theorem \ref{thm:slowEntRankone1}) proved that the rank-one system has growth rate at most $n^2$. It is worth to compare  Theorem \ref{thm:entropyConvergenceRateRankOne} with Theorem \ref{thm:slowEntRankone1} as they show that complexity of rank-one system has the same bounds.
\erm

\brm
By considering information function $I_{\alpha}(x)=-\log_2\mu(\alpha(x))$, we can formalize pointwise version results parallel to Theorem \ref{thm:entropyConvergenceRateAperiodic} and Corollary \ref{cor:entropyConvergenceRateComErgo}.
\erm




\section{Open questions}\label{sec:openQuestions}
In this section we propose several open questions related to slow entropy type invariants.

\subsection{Slow entropy and spectral properties}
Recall that a measure-preserving transformation $T$ on $(X,\mathscr{B},\mu)$ is said to have \emph{countable Lebesgue spectrum} if there exists an infinite orthonormal set $\{1,\phi_0,\ldots\}\subset L^2(X)$ such that $\phi_{ik}=U_T^k(\phi_i)\in1^{\bot}$ are all pairwise distinct and orthogonal, where $U_T$ is the Koopman operator induced from $T$. It is interesting if there is any relation between  between slow entropy and Lebesgue spectrum:
\begin{question}
Does countable Lebesgue spectrum imply any below restrictions on the slow entropy?
\end{question}

Moreover, one can ask the following question for rigid systems (recall that rigidity implies singular spectrum):
\begin{question}
Is it possible to have the upper slow entropy for a rigid transformation positive with respect to $a_n(t)=n^t$?
\end{question}

\subsection{Slow entropy and finite rank systems}
Recall that a system is a finite rank system if instead of approximating it with one sequence of towers (which $\epsilon$-refines every partition) one is allowed finitely many such towers\footnote{For more precise definition, see \cite{Ferenczirank}.} We have the following question:

\begin{question}\label{question:finiteRank}
What about the finite rank systems' slow entropy?
\end{question}
It is worth to point that in general rank one situation Ferenczi \cite{Ferenczi} and Kanigowski \cite{Kanigowski} obtain estimates of upper bounds of lower slow entropy, see Theorem \ref{thm:slowEntRankone1} and Theorem \ref{thm:rankoneSlowEntropy} for more details.

Recall that a very natural class of finite rank systems is given by interval exchange transformations (IET's). More precisely, suppose $n>0$, $\pi$ is a permutation defined on $\{1,\ldots,n\}$ and $\lambda=(\lambda_1,\ldots,\lambda_n)$, where $\sum_{i=1}^n\lambda_i=1$ and $\lambda_i>0$ for $1\leq i\leq n$. Then a map $T_{\pi,\lambda}:[0,1]\to[0,1]$ is an $(\pi,\lambda)$-interval exchange transformation if for $x\in[a_i,a_i+\lambda_i)$, we have  $$T_{\pi,\lambda}(x)=x-a_i+a'_i,$$ where $a_i=\sum_{1\leq j\leq i}\lambda_j$, $a'_i=\sum_{1\leq j\leq\pi(i)}\lambda_{\pi^{-1}(j)}$ for $1\leq i\leq n$.

It is interesting to study growth rate for IET's:

\begin{question}
For an interval exchange transformation, is that possible its upper (lower) slow entropy is $1$ in the scale $a_n(t)=n^t$?\end{question}

However, if we restricted ourselves in the rank one situation, by adding some additional ergodic conditions, we may ask the following simple version of Question \ref{question:finiteRank}\footnote{See conjecture in \cite{Ferenczi} for more explanations.}:
\begin{question}
What about the slow entropy of Ornstein \cite{Orn} mixing rank one system? While it must satisfy Theorem \ref{thm:rankoneSlowEntropy}, do they have an infinite upper slow entropy when $a_n(t)=n^{kt}$ for any $k$?
\end{question}

\subsection{Realization of slow entropy}
One of important combinatorial tools in ergodic theory is fast approximation (See \cite{KatokCombinatorial} for more details.). Thus it is reasonable to consider the realization of ergodic systems with certain given slow entropy by fast approximation:
\begin{question}
Is it possible to have a system constructed by fast periodic approximation with any subexponential slow entropy (upper or lower)?
\end{question}

\subsection{Other related questions}
Recall that slow entropy is designed to measure the complexity of zero entropy ergodic systems, thus the estimates of certain representative systems slow entropy would be quite interesting and meaningful:
\begin{question}
Can we get estimates of the slow entropy for billiards\footnote{For detailed discussion about billiards in polygons, we refer to \cite{ForniMatheus}.} in polygons and polyhedra?
\end{question}

\begin{question}
Can we get estimates of the slow entropy of Gaussian\footnote{For detailed discussion of Gaussian system, we refer to \cite{LemParThou}.} systems?
\end{question}

Finally we recall the following (well-known) question on a relation between classical entropy and sequence entropy:
\begin{question}
Does there exist a zero-entropy system whose sequence entropy is positive with respect to the sequence $\{n^2\}$. One can ask the same question for the sequence $\{p(n)\}$ where $p$ is a polynomial.
\end{question}

\end{document}